\input amstex
\documentstyle{amsppt}
\magnification=1200

\topmatter
\title
Incidence coefficients in the
Novikov Complex for Morse forms: rationality and
exponential growth properties.
\endtitle
\leftheadtext{A.Pajitnov}
\rightheadtext{Incidence coefficients in the Novikov complex}
\author
A.Pajitnov
\endauthor
\address
Universit\'e de Nantes, D\'epartement de Math\'ematiques,
2 rue de la Houssini\`ere, 44072 Nantes Cedex
\endaddress
\email
pajitnov\@univ-rennes1.fr
\endemail

\define\a{\alpha}
\redefine\b{\beta}
\redefine\d{\delta}
\define\k{\varkappa}
\redefine\t{\theta}
\redefine\l{\lambda}
\define\m{\mu}
\define\om{\omega}

\define\s{\sigma}
\define\g{\gamma}
\define\e{\epsilon}

\redefine\o{\omega}

\define\ffmin{f^{-1}}
\define\vvol{\text{\rm vol}}

\define\odin{\bold 1}
\define\eg{~of~exponential~growth~}

\define\Mat{\text\rm Mat}
\define\Tub{\text{\rm Tub}}

\define\QQQ{\bold Q}
\define\LL{\Cal L}

\redefine\D{\Delta}

\redefine\L{\Lambda}

\define\G{\Gamma}

\redefine\O{\Omega}

\define\vxi{v\langle\vec\xi\rangle}

\define\ma{~manifold~}
\define\gr{gradient~}
\define\tr{~trajectory~}

\define\grs{gradients~}
\define\trs{~trajectories~}

\define\Mf{~Morse~function~}
\define\co{~cobordism~}
\define\sma{~submanifold~}
\define\ho{~homomorphism~}

\define\th{~therefore~}
\define\Th{~Therefore~}

\define\ic#1{~integral~{#1}-cone}
\define\sq#1{[{#1}]}
\define\gsq#1{\ZZZ{\sq{#1}}}
\define\noconf{~there~is~no~possibility~of~confusion~}
\define\acc#1{\{#1\}}
\define\wi{\widetilde}
\define\bbar#1{\overline{\overline #1}}

\redefine\o{\omega}

\define\zd{\ZZZ[Z]}

\define\zzz{\ZZZ[\ZZZ^m]}

        \TagsOnLeft
\TagsAsMath

\define\h{\widehat}
\define\qu{\QQ^{-1}}
\define\sm{\setminus}

\define\card{\text{\rm card}}

\redefine\leq{\leqslant}
\redefine\geq{\geqslant}

\redefine\AA{{\Cal A}}
\define\PP{{\Cal P}}
\define\QQ{{\Cal Q}}

\define\VV{{\Cal V}}

\define\UU{{\Cal U}}

\define\HH{{\Cal H}}
\define\GG{{\Cal G}}

\define\gt{{\Cal G}t}

\define\aand{\quad\text{and}\quad}
\define\wwhere{\quad\text{where}\quad}
\define\ffor{\quad\text{for}\quad}

\define\iif{\quad\text{if}\quad}

\define\RRR{{\bold R}}
\define\ZZZ{{\bold Z}}
\define\NNN{{\bold N}}

\define\Wkr{W^{\circ}}

\define\Ker{\text{\rm Ker }}
\define\ind{\text{\rm ind}}
\define\rk{\text{\rm rk }}
\redefine\Im{\text{\rm Im }}
\define\supp{\text{\rm supp }}
\define\Int{\text{\rm Int }}

\define\smo{C^{\infty}}

\define\fpr#1#2{{#1}^{-1}({#2})}
\define\sdvg#1#2#3{\widehat{#1}_{[{#2},{#3}]}}
\define\disc#1#2#3{B^{({#1})}_{#2}({#3})}
\define\Disc#1#2#3{D^{({#1})}_{#2}({#3})}
\define\desc#1#2#3{B_{#1}(\leq{#2},{#3})}
\define\Desc#1#2#3{D_{#1}(\leq{#2},{#3})}
\define\komp#1#2#3{{\bold
K}({#1})^{({#2})}({#3})}
\define\Komp#1#2#3{\big({\bold
K}({#1})\big)^{({#2})}({#3})}

\define\ran{\{(A_\lambda , B_\lambda)\}_{\lambda\in\Lambda}}
\define\rans{\{(A_\sigma , B_\sigma)\}_{\sigma\in\Sigma}}

\define\rran{\{(A_\lambda^{(s)} , B_\lambda^{(s)})\}_{\lambda\in\Lambda,
0\leq s\leq n}}
\define\rrans{\{(A_\sigma^{(s)} , B_\sigma^{(s)})\}_{\sigma\in\Sigma,
0\leq s\leq n}}

\define\fmin{F^{-1}}

\define\chart{\Phi_p:U_p\to B^n(0,r_p)}
\define\atlas{\{\Phi_p:U_p\to B^n(0,r_p)\}_{p\in S(f)}}
\define\alas{\{\Phi_p:U_p\to B^n(0,r_p)\}_{p\in S(\o)}}

\define\stind#1#2#3{{#1}^{\displaystyle\rightsquigarrow}_
{[{#2},{#3}]}}

\define\indl#1{{\scriptstyle{\text{\rm ind}\leqslant {#1}~}}}
\define\inde#1{{\scriptstyle{\text{\rm ind}      =   {#1}~}}}

\define\id{\text{id}}

\define\st#1{\overset\rightsquigarrow\to{#1}}
\define\bst#1{\overset{\displaystyle\rightsquigarrow}\to{\boldkey{#1}}}

\define\stexp#1{{#1}^{\rightsquigarrow}}
\define\bstexp#1{{#1}^{\displaystyle\rightsquigarrow}}

\define\bstind#1#2#3{{\boldkey{#1}}^{\displaystyle\rightsquigarrow}_
{[{#2},{#3}]}}
\define\bminstind#1#2#3{\stind{({\boldkey{-}\boldkey{#1}})}{#2}{#3}}

\define\kr#1{{#1}^{\circ}}

\abstract
In this paper we continue the study of generic properties 
of the Novikov complex, began in  the work
"The incidence coefficients in the Novikov complex
are generically rational functions" ( dg-ga/9603006).

For a Morse map $f:M\to S^1$ there is a refined version
of Novikov complex, defined over the Novikov completion
of the fundamental group ring. We prove that for a 
$C^0$ generic $f$-gradient the corresponding
 incidence coefficients belong to the image in the Novikov
ring of a (non commutative) localization of
the fundamental group ring.

The Novikov construction generalizes also to the
case of
 Morse 1-forms. In this case the corresponding
incidence coefiicients belong to
 the suitable completion of the ring of
integral Laurent polynomials of several variables.
We prove that for a given Morse form
$\o$ and a $C^0$ generic $\o$-gradient 
 these incidence coefficients
are rational functions.

The incidence coefficients in the Novikov
complex
are obtained by counting the
algebraic number of the trajectories of the gradient,
 joining
the zeros of the Morse form.
There is V.I.Arnold's version of the exponential growth
conjecture, which concerns the total number
of trajectories. Namely,  $\o$ be a Morse form on 
a closed manifold, $v$ be an $\o$-gradient,
$p:M'\to M$ be a free abelian covering for which
$p^*\o =df$ with $f:M'\to\RRR$. Let $x$ be a critical point
of $f$ of index $k$ and $c$ be a real number.
 The conjecture says that the number of $v$-trajectories
 joining $x$ to the critical
points $y$ of $f$ of index $k-1$ with $f(y)\geq c$ grows at most
exponentially with $-c$.
 We confirm it  for any given
Morse form and a $C^0$ dense set of its gradients.

We give an example of explicit 
computation of the Novikov complex.  
\endabstract
\endtopmatter

\document

\head{ Introduction} \endhead
\subhead{A.  Morse-Novikov theory}\endsubhead
 The classical Morse-Thom-Smale construction
associates to a Morse function $g:M\to\RRR$
on a closed manifold a free chain complex
$C_*(g)$ where the number $m(C_p(g))$ of free generators
of $C_p(g)$ equals the number of the critical points of $g$ of
index $p$ for each $p$. The boundary operator in this
complex is defined in a geometric way, counting the trajectories
of a gradient of $g$, joining critical points of
$g$ (see [4], [7], [10], [11], [12]).

In the early 80s S.P.Novikov generalized this construction
 to the case of maps $f:M\to S^1$ ( see [5]).
The corresponding analog of Morse complex
is a free chain complex $C_*(f)$ over
$\ZZZ [[t]][t^{-1}]$. Its number of free generators
equals the number of critical points of $f$ of index $p$,
and the homology of $C_*(f)$ equals to the completed
homology of the cyclic covering.

Fix some $k$. The boundary operator
 $\partial :C_k(f)\to C_{k-1}(f)$
is represented by a matrix, which entries are in
the ring of Laurent power series. That is
$\partial_{ij}=\sum_{n=-N}^{\infty} a_nt^n,\wwhere
\mathbreak
 a_n\in\ZZZ$.

Since the beginning S.P.Novikov conjectured that the
power series
$\partial_{ij}$
 have some nice analytic
properties. In particular he conjectured that
\vskip0.1in
\it
Generically the coefficients $a_n$ of
$\partial_{ij}=\sum_{n=-N}^{\infty} a_nt^n$
grow at most exponentially with $n$.
\rm
\vskip0.1in
In [9] we have proved that for a $C^0$ generic $f$-\gr
the incidence coefficients above are actually
rational functions. To recall the statement of the
Main Theorem of [9] let
 $M$ be a closed connected manifold and $f:M\to S^1$
a Morse map, non homotopic to zero. Denote the set of
critical points of $f$ by $S(f)$.
The set of  $f$-gradients
of the class  $\smo$,
 satisfying the transversality
assumption (see \S1 for terminology),
 will be denoted by ${\Cal G} t(f)$.
By Kupka-Smale theorem it is residual in the set of
all the $\smo$ gradients. Choose
$v\in\gt (f)$. Denote by $\bar M @> {\Cal P} >>  M$
the connected infinite cyclic covering for which
$f\circ \Cal P$ is homotopic to zero.
Choose a lifting  $F:\bar M\to\RRR$ of
$f\circ \Cal P$ and let $t$ be the generator
of the structure group of $\Cal P$
such that $F(xt)<F (x)$.
The $t$-invariant lifting of $v$ to $\bar M$
will be denoted by the same letter $v$.
For every critical point $x$ of $f$ choose a
 lifiting $\bar x$
of $x$ to $\bar M$. Choose orientations
of stable manifolds of critical points.
Then for every $x,y\in S(f),~\ind x =\ind y +1$
and every $k\in\ZZZ$ the incidence coefficient
$n_k(x,y;v)$
is defined
(as the algebraic number of $(-v)$-trajectories
joining $\bar x$ to $\bar y t^k$).

\proclaim{ Theorem [9,p.2]}
In the set $\gt (f)$ there is a subset $\gt _0(f)$ with
the following properties:
\roster\item
$\gt_0 (f)$ is open and dense in $\gt (f)$ with respect to
$C^0$ topology.
\item If $v\in\gt _0(f),~x,y\in S(f)\aand \ind x=\ind y+1$,
then
$\sum_{k\in\ZZZ} n_k(x,y;v)t^k$ is a rational
function of $t$ of the form
$\frac{P(t)}{t^m Q(t)}, \wwhere
 P(t) \aand Q(t)$
are polynomials with integral
coefficients,
$m\in\NNN$,~and~$Q(0)=1$.
\item Let $v\in\gt _0(f)$.
Let $U$ be a neighborhood of $S(f)$.
Then for every $w\in\gt _0(f)$ such that
$w=v$ in $U$ and $w$
is sufficiently close to $v$ in $C^0$ topology
we have:
$n_k(x,y;v)=n_k(x,y;w)$
for every $x,y\in S(f),~k\in\ZZZ$.
\endroster
\endproclaim

In the present paper we develop the methods of [9]
and apply them to the incidence coefficients with values in the
Novikov completion of the fundamental group ring
$\ZZZ\pi_1M$ (\S 1). In  \S 2 we consider the case
of arbitrary Morse forms and free abelian coverings.
We also give an example of explicit computation of the
Novikov complex (\S 3). In \S 4 we prove the 
V.I.Arnold conjecture, concerning the
total number of $(-v)$-trajectories, joining the critical points
of adjacent indices.

\subhead
B.  Morse forms and Novikov rings
\endsubhead
To give the statement of our results, we recall some algebraic
and Morse-theoretic definitions.

Let $G$ be  a group and $\xi:G\to\RRR$ be a group \ho.
We denote by $(\ZZZ G)\sphat~\sphat$ the abelian group
of all formal linear combinations
$\sum\limits_{g\in G}n_gg$ (infinite in general).
Recall that the Novikov ring $\ZZZ G_\xi^-$ is the
ring of such $\l\in (\ZZZ G)\sphat~\sphat,~
\l=\sum_{g\in G}n_gg$,
that for every $c\in\RRR$ the set $\supp \l\cap \xi^{-1}([c,\infty[)$
is finite.

Let $\o$ be a closed 1-form on a \ma $M$.
The deRham cohomology class of $\o$ will be denoted
by $[\o]$ and
 the corresponding \ho $\pi_1M\to\RRR$
will be denoted by $\acc\o$.
 We say that $\o$ is a
\it Morse
form  \rm , if locally it is the differential of a Morse function.
If $f:M\to S^1$ is a Morse map, then its differential
is a Morse form, which cohomology class is in $H^1(M,\ZZZ)$.
A Morse form $\o$ is proportional to
a differential of a Morse map $M\to S^1$ if and only if
$\exists\l\in\RRR : \l[\o]\in H^1(M,\ZZZ)$.

The terminology of \S 1.A of [9] (which concerns Morse functions)
is extended in an obvious way to the case of Morse forms, and
we shall make free use of it. In particular we shall assume the
notion of \it $\o$-\gr\rm.
The set of all $\o$-gradients, satisfying the transversality assumption
will be denoted by $\GG t(\o)$.

Let $\o$ be a Morse form on a closed connected manifold $M$
and let $v\in\GG t(\o)$.
Let $x,y\in S(\o), \ind x=\ind y+1$. Choose some liftings
$\wi x,\wi y$ of $x,y$ to $\wi M$ and the orientations of the
stable manifolds of $x$ and of $y$. Then the incidence
coefficient
$\wi n(\wi x,\wi y;v)\in\ZZZ(\pi_1M)_\xi^-$
is defined (for the case of  Morse maps $M\to S^1$
see the precise definition of $\wi n(\wi x,\wi y;v)$
in [7];  the general case follows by the approximation
procedure;
see Lemma 2.6 and the discussion before it).

\subhead
C. Statement of the results
\endsubhead
\vskip0.2in
\it 1. Morse maps $M\to S^1$
\rm \vskip0.2in

Let $\xi:G\to\ZZZ$ be a group epimorphism.
 Denote $\Ker\xi$
by $H$. For $n\in\ZZZ$ denote $\xi^{-1}(n)$ by $G_{(n)}$
and
$\{x\in\ZZZ G\mid\supp x\subset G_{(n)}\}$
by $\ZZZ G_{(n)}$.
Denote $\xi^{-1}(]-\infty,-1])$ by $G_-$ and
$\{x\in\ZZZ G\mid\supp x\subset G_-\}$
by $\ZZZ G_-$. Choose $\t\in\ZZZ G_{(-1)}$.
It is easy to see that $(\ZZZ G)^-_\xi$ is identified with the
ring of power series of the form
$\{a_{-n}\t^n+...+a_1\t+... \mid a_i\in\ZZZ H\}.$

Set
$\Sigma_n=\{\odin+A\mid A\in\Mat_{n\times n}(\ZZZ G_{(-1)})\}$.
   Set  $\Sigma=\bigcup\limits_{n\geq 1}\Sigma_n$.

There is the corresponding localization ring
$\ZZZ G_\Sigma$ (see [3, p.255]). Every matrix in
$\Sigma_n$ is invertible in $\Mat_{n\times n}(\ZZZ G_\xi^-)$,
the inverse of $\odin+A$ being given by
$\sum_{n=0}^\infty (-1)^n A^n$, \th the localization
map $\l:\ZZZ G\to\ZZZ G_\Sigma$ is injective
and the inclusion $i:\ZZZ G\hookrightarrow \ZZZ G_\xi^-$
factors through a ring \ho
$\ell:\ZZZ G_\Sigma\to\zzz _\xi^-$.

Let $M$ be a connected closed manifold and $f:M\to S^1$
be a Morse map, nonhomotopic to zero.
Denote by $\xi$ the induced \ho $\pi_1M\to\ZZZ$.
Denote by $p:\wi M\to M$ the universal covering of $M$.

\proclaim{Theorem A}
In the set $\gt (f)$ there is a subset $\gt _1(f)$ with
the following properties:
\roster\item
$\gt_1 (f)$ is open and dense in $\gt (f)$ with respect to
$C^0$ topology.
\item If $v\in\gt _1(f)$
then
for every $x,y\in S(f)$ with $\ind x=\ind y+1$ we have
$\wi n(\wi x,\wi y;v)\in \Im\ell$.
\item Let $v\in\gt _1(f)$.
Let $U$ be a neighborhood of $S(f)$.
Then for every $w\in\gt _1(f)$ such that
$w=v$ in $U$ and $w$
is sufficiently close to $v$ in $C^0$ topology
we have:
$\wi n(\wi x,\wi y;v)=\wi n(\wi x,\wi y;w)$
for every $x,y\in S(f)$.
\endroster
\endproclaim

\vskip0.2in\it
2. Morse forms with arbitrary cohomology classes \rm
\vskip0.2in

At present we can prove the analog of the Theorem A
in the case of arbitrary Morse forms only
for the incidence coefficients associated with
free abelian coverings.

Let $\o$ be a Morse form on a 
closed connected
\ma $M$. If $\phi:\h M\to M$ is any
regular covering with the structure group $G$, such that
$\phi^*([\o])=0$, then the \ho
$\acc\o:\pi_1M\to\RRR$ factors as
$\pi_1M\to G \to\RRR$ and it is not difficult to see that
the incidence coefficients $\h n(\h x,\h y;v)$ are defined for
every $v\in\GG t(\o)$
(here we suppose that $\ind x=\ind y+1$, and that
for every $p\in S(\om)$ a lifting
$\hat p$ of $p$ to $\h M$ and an orientation
of the stable \ma of $p$ are chosen).
 In particular it is the case for
the maximal free abelian covering
$\bbar M@>\PP >> M$ with the structure group
$H_1(M,\ZZZ)/\text{\rm Tors}\approx \ZZZ^m$.
By abuse of notation we shall denote the corresponding
\ho $\ZZZ^m\to\RRR$ by the same symbol 
as the de Rham cohomology class 
$[\o]$ of $\o$. Assume that $[\o]\not=0$.
Set
$S_{[\o]}=\{P\in\ZZZ[\ZZZ^m] ~\vert~ P={\bold 1}+Q: \supp Q\subset
[\o]^{-1}(]-\infty,0[)\}$.

\proclaim{Theorem B}
There is a subset
$\GG t_1(\o)\subset\GG t(\o)$ with
the following properties:
\roster\item
$\GG t_1(\o)$ is open and dense in $\GG t(\o)$ with respect to
$C^0$ topology.
\item
For every $v\in\GG t_1(\o)$ and every $x,y\in S(\o)$
with $\ind x=\ind y+1$ we have:
$\bbar n(\bbar x,\bbar y;v)\in S_{[\o]}^{-1}\zzz$.
\item Let $v\in\gt _1(\o)$.
Let $U$ be a neighborhood of $S(\o)$.
Then for every $w\in\GG t_1(\o)$ such that
$w=v$ in $U$ and $w$
is sufficiently close to $v$ in $C^0$ topology
we have:
$\bbar n(\bbar x,\bbar y;v)=\bbar n(\bbar x,\bbar y;w)$
for every $x,y\in S(\o),~k\in\ZZZ$.
\endroster
\endproclaim

\vskip0.2in
\it 3. Exponential growth estimates \rm
\vskip0.2in

Let $G$ be a group.
For an element $a=\sum n_gg\in\ZZZ G$ we denote by $\Vert a\Vert$
the sum $\sum\vert n_g\vert$.

Let $\xi:G\to\RRR$ be
a \ho. For 
$\l=\sum\limits_{g}n_gg\in\ZZZ G_\xi^-$ and $c\in\RRR$ we denote by
$\l[c]$ the element
$\sum\limits_{\xi(g)\geq c}n_gg$ of $\ZZZ G$ and we set
$N_c(\l)=\Vert\l[c]\Vert$. We shall say that $\l$ is of exponential growth
if there are $A,B\geq0$ such that for every $c\in\RRR$
we have $N_c(\l)\leq Ae^{-cB}$. It is easy to prove
that the elements of exponential growth form a subring of
$\ZZZ G_\xi^-$, which contains $\ZZZ G$.

\proclaim{Theorem C}
Let $v$ be an $\o$-gradient,
belonging to $\GG t_1(\o)$.
Let $x,y\in S(\o),\ind x=\ind y+1$. Then
$    \wi n(\wi x,\wi y;v) \in\ZZZ[\pi_1M]^-_{\{\o\}}$ is \eg.
\endproclaim

\vskip0.2in\it
4. An example\rm
\vskip0.2in

In the subsection 3 we construct a three-manifold $M$,
 a Morse map $f:M\to S^1$ and an $f$-gradient $v$
such that
$n_0(\bar x,\bar y;v)=0$ and for $k\geq 0$ we have
 $n_{k+1}(\bar x,\bar y;v) =
-\frac 4{\sqrt 5}\cdot \bigg( \bigg(\frac {3+\sqrt 5}2\bigg)^k -
 \bigg(\frac {3-\sqrt 5}2\bigg)^k\bigg).$

\vskip0.2in\it
5. Exponential estimates of absolute
number of trajectories: Morse maps $M\to S^1$
\rm
\vskip0.2in

We assume here the terminology of Subsection A. The set
of all
$f$-gradients of class $C^{\infty}$ will be denoted by $\GG (f)$.
Recall from [8,\S 2B] that an $f$-\gr $v$ is called
\it good \rm
if for every $p,q\in S(f)$ we have
$$\bigg(\ind p\leq \ind q+1\bigg)
\Rightarrow
\bigg(D(p,v)\pitchfork D(q,-v)\bigg)$$

The set of all good $f$-gradients will be denoted by
$\GG d(f)$. For $v\in \GG (f)$ we denote by the same letter $v$
the $t$-invariant lifting of $v$ to $\bar M$.
Choose a lifting $F:\bar M\to\RRR$ of $ f$ to $\bar M$.
It is easy to prove that for
$p,q\in S(F),\ind p=\ind q+1$ and for for $v\in \GG d(f)$ the set of
$(-v)$-trajectories, joining $p$ to $q$ is finite.
The liftings of critical points of $f$ to $\bar M$
being chosen, denote by $N_k(x,y;v)$ the number of
$(-v)$-trajectories joining $\bar x$ to $\bar yt^k$
(where $\ind x=\ind y+1$).

\proclaim{Theorem D} In the set $\GG (f)$ there is a subset
$ \GG_0(f)$ with the following properties:
\roster\item $\GG _0(f)$ is $C^0$ dense in $\GG (f)$ and
$\GG_0(f)\subset \GG d(f)$.
\item Let $v\in \GG_0(f)$.
Then there are constants $C,D>0$ such that for every
$x,y\in S(f)$ with $\ind x=\ind y+1$ and for every $k\in\ZZZ$ we have
$N_k(x,y;v)\leq C\cdot D^k$.    \endroster
\endproclaim

\vskip0.2in\it
6. Exponential estimates of absolute
number of trajectories: Morse forms
\rm
\vskip0.2in

Let $\om$ be a Morse form on a closed connected manifold $M$.
Let $\PP:\bbar M\to M$ be the maximal free abelian covering of
$M$ with structure group $\ZZZ^m$; we identify
the cohomology class $[\om]$ of $\om$ with the corresponding \ho
$\ZZZ^m\to\RRR$.
We denote by $\GG (\om)$ the set of all $\om$-gradients
of class $C^{\infty}$ and by $\GG d(\om)$
the set of all good $\om$-\grs of class $C^{\infty}$.
Let $v\in\GG d(\om)$.
For every zero $x$ of $\om$ choose a
 lifiting $\bbar x$
of $x$ to $\bbar M$ and an orientation
of the stable manifold of $x$.
Then for every $g\in\ZZZ^m$ and every $x,y\in S(\om)$ with $\ind x=\ind y+1$
the set of $(-v)$-\trs joining $\bbar x$ to $\bbar yg$ is finite
and we denote its cardinality
by $N(\bbar x,\bbar y,g;v)$. For $c\in\RRR$ we denote by
$N_{\geq c}(\bbar x,\bbar y;v)$ the sum
$\sum\limits_{g: [\om](g)\geq c} N(\bbar x,\bbar y;g;v)$.

\proclaim{Theorem E}
In the set $\GG (\om)$ there is a subset $\GG_0(\om)$ with the
following properties:
\roster\item $\GG_0(\om)$ is dense in $\GG (\om)$ with respect
to $C^0$ topology; $\GG_0(\om)\subset \GG d(\om)$.
\item Let $v\in \GG_0(\om)$. There are constants $C,D>0$ such that for every
$x,y\in S(\o)$ with $\ind x=\ind y+1$ and every
$c\in\RRR$ we have
$N_{\geq \l}(x,y;v)\leq C\cdot D^{-\l}$.  \endroster
\endproclaim

\newpage

\newpage

\head
\S 1.  Morse maps $M\to S^1$
\endhead
\subhead
A. Algebraic preliminaries
\endsubhead

We accept here the terminology of Subsection C.1 of Introduction.
Further, we say that an element $\k\in\ZZZ(\pi_1M)^-_\xi$
is of type $(\Cal L)$, if
\vskip0.2in

$(\Cal L)$ There are $r,q\in G$, a natural number $m$,
an $m\times m$-matrix $\AA=(a_{ij})_{1\leq i,j\leq m}$ where
$a_{ij}\in\ZZZ G_{(-1)}$ and vectors
$(X_i)_{1\leq i\leq m},
 (Y_i)_{1\leq i\leq m}, X_i,Y_i\in\ZZZ H$,
such that
$$\k=r\bigg(\sum\limits_{s\geq 0}\sum \Sb  \\ 1\leq i\leq m \\
1\leq j\leq m\endSb Y_ia_{ij}^{(s)}X_j\bigg)q$$
\vskip0.2in
where $a_{ij}^{(s)}$ are the entries of ${\Cal A}^s$.
The next lemma follows from the definition of the \ho $\ell$.

\proclaim{Lemma   1.1}
The elements of type $(\Cal L)$ are contained in $\Im \ell$.
\quad$\square$\endproclaim

\proclaim{Proposition 1.2}
The elements of type $(\Cal L)$
are \eg.
\endproclaim
\demo{Proof}
It suffices to prove that every matrix entry
of the matrix series
$u=\sum_{s\geq 0}A^s$, where $A=(a_{ij})$ and
$a_{ij}\in\ZZZ G_{(-1)}$, is \eg.
For an $(m\times m)$-matrix $B=(b_{ij})$ we denote
by $\Vert B\Vert$ the number $\max\limits_{i,j}\Vert b_{ij}\Vert$.
It is easy to check that
$\Vert BC\Vert\leq\Vert B\Vert\cdot\Vert C\Vert\cdot m$.
Let $A$ be $(m\times m)$-matrix. Then
$\Vert A^s\Vert\leq\Vert A\Vert^s\cdot m^{s-1}\leq
\Vert A\Vert^s\cdot m^{s}$.
Let $1\leq i,j\leq m$. Write $u_{ij}=\sum n_gg$.
Then for $k\geq 0$ we have
$\sum\limits_{\xi(g)=-k}\vert n_g\vert\leq (m\cdot\Vert A\Vert)^k$.
If $m\Vert A\Vert\leq 1$ this gives
$\sum\limits_{\xi(g)\geq -k}\vert n_g\vert\leq k+1\leq e^k$.
If $m\Vert A\Vert> 1$ we have
$\sum\limits_{\xi(g)\geq -k}\vert n_g\vert\leq
\frac {(m\Vert A\Vert)^k-1}{(m\Vert A\Vert)-1}
<D(m\Vert A\Vert)^k$.

\Th in any case there are
$A,B>0$ such that for $k>0, k\in\ZZZ$ we have
$\sum\limits_{\xi(g)\geq k}\vert n_g\vert
\leq A\cdot B^{-k}$. For $k\geq 0$ it is true obviously and
this implies
that $u$ is \eg.
\quad$\square$\enddemo

\subhead
B. Statement of Theorem 1.3
\endsubhead

 Theorem A follows immediately from the next theorem.

\proclaim{Theorem 1.3}
In the set $\gt (f)$ there is a subset $\gt _1(f)$ with
the following properties:
\roster\item
$\gt_1 (f)$ is open and dense in $\gt (f)$ with respect to
$C^0$ topology.
\item If $v\in\gt _1(f)$
then
for every $x,y\in S(f)$ with $\ind x=\ind y+1$ we have:
$\wi n(\wi x,\wi y;v)$ satisfies $(\LL)$.
\item Let $v\in\gt _1(f)$.
Let $U$ be a neighborhood of $S(f)$.
Then for every $w\in\gt _1(f)$ such that
$w=v$ in $U$ and $w$
is sufficiently close to $v$ in $C^0$ topology
we have:
$\wi n(\wi x,\wi y;v)=\wi n(\wi x,\wi y;w)$
for every $x,y\in S(f)$ such that $\ind x=\ind y+1$.
\endroster
\endproclaim

\subhead
C. Generalities on intersection indices
\endsubhead
Let $M$ be a manifold without boundary, $\QQ : \h M\to M$
be a regular covering (not necessarily connected)
with structure group $H$. We say, that a submanifold $N$
of $M$ is \it lifted-oriented \rm
(resp. \it lifted-cooriented \rm )
if a lifting $\hat i :N\hookrightarrow \h M$
of the inclusion map $i:N\hookrightarrow M$ is fixed,
and $N$ is oriented (resp. cooriented). We shall
denote $\h i (N)$ by $\h N$.

Let $X\subset M$ and let $N$ be a lifted-oriented
submanifold of $M$ such that $N\sm \Int X$ is compact.
Then $\h N\sm\Int\qu(X)$ is compact,
and the homology class $[\h N]_{\h M,\qu (X)}\in
H_n(\h M,\qu(X))$ is defined, where $n=\dim N$
(see [9,\S 4]). We shall denote it by $[\h N]_{M,X}$,
or simply by $[\h N]$ if there is no posibility of confusion.
Let $L$ is a compact lifted-cooriented submanifold without
boundary of $M$. Then there is the coorientation class
$]\h L[\in
\mathbreak
H^{m-l}(\h M,\h M\sm\h L)$, where
$l=\dim L$, $m=\dim M$.

Assume that $X\cap L=\emptyset, N\pitchfork L$
and $n+l=m$. Then $\h N\pitchfork \h L$, and
$\h N\cap \h L$ is finite, and the intersection index
$\h N~\sharp~ \h L\in\ZZZ$
is defined. Denote by $j$ the inclusion
$(\h M,\QQ^{-1}(X))\hookrightarrow (\h M,\h M\sm\h L)$.
The next lemma is standard.
\proclaim{Lemma 1.4}
$\h N~\sharp~\h L=j^*(]\h L[)([\h N])$.\quad$\square$
\endproclaim

\subhead
D. Ranging systems
\endsubhead
Let $f:W\to[a,b]$ be a Morse function on a compact
riemannian cobordism, $f^{-1}(b)=V_1, f^{-1}(a)=V_0$,
$v$ be an $f$-gradient. Let $\QQ:\h W\to W$ be a regular covering with
a structure group $H$. The lifting of $v$ to $\h W$ will be denoted by
$\h v$. If $x\in W$ and $\g(x,t;v)$ is defined on $[0,a]$, and $\h x\in\QQ^{-1}(x)$,
then the lifting to $\h W$ of $\g(x,\cdot;v)$, starting at $\h x$
is the $\h v$-trajectory $\g(x,\cdot;\h v)$. It
 is easy to define with the help of this lifiting
procedure a diffeomorphism
$\stexp {\h v} :\qu(V_1\setminus K(-v))\to \qu(V_0\setminus K(v))$.
For $X\subset V_1$ we denote by abuse of notation
$\stexp {\h v} (X\setminus K(-v))$ by $\stexp {\h v}(X)$.

Let $N$ be a oriented-lifted submanifold of $V_1$. Then it is easy to see
that
$\st v(N)$ is a oriented-lifted submanifold of $V_0$.

Now let $\ran$ be  a ranging system for $(f,v)$
 (see [9,\S 4 Subsection B] for
definitions). Let $N$ be an oriented-lifted submanifold
of $V_1\setminus B_b$ such that $N\setminus \Int A_b$ is compact.
Then Prop. 4.6 of [9] implies that $\st v (N)$ is an oriented-lifted submanifold of
$V_0\setminus B_a$ such that $\st v(N)\setminus \Int A_a$ is compact.

The following proposition is a generalization of the Proposition
4.7 of [9], which can be considered as a particular case ($H=\{1\}$)
of the Proposition 1.5. The proof of 1.5 is carried out along the lines
of [9, \S 4, Subsection B]. We recommend to the reader to consult 
[9,\S 4, Subsection B] for the basic definitions (such as the definition
of ranging system, cited above ). We present below the
main steps of proof.

\proclaim{Proposition 1.5}
There is a homomorphism
$\h H(v):H_*(\qu(V_1\sm B_b),\qu (A_b))\to
H_*(\qu(V_0\sm B_a),\qu (A_a))$
of right $\ZZZ H$-modules, such that:

\roster\item If $N$ is an oriented-lifted submanifold of
$V_1\sm B_b$, such that $N\sm\Int A_b$ is compact, then
$\h H(v)([\h N])=
\big[(\st v (N))\sphat~\big]. $
\item There is an $\e>0$ such that for every $f$-gradient $w$
with $\Vert w-v\Vert<\e$ we have $\h H(v)=\h H(w)$.
\endroster
\endproclaim

\it Proof. \rm\quad
 An easy induction argument shows that it is sufficient
to prove the proposition in the case $\card \Lambda=1$.
Let $S(f)=S1(f)\sqcup S2(f)$, where
for every $p\in S1(f)$, resp. $p\in S2(f)$
the i), resp. ii) of (RS2)
of Definition 4.3 in [9,p.23]
 holds.
Pick Morse functions $\phi_1, \phi_2:W\to[a,b]$, adjusted
to $(f,v)$, such that there are regular values
$\mu_1$ of $\phi_1$, $\mu_2$ of $\phi_2$
satisfying:
(1) for every $p\in S1(f)$ we have: $\phi_1(p)<\mu_1$
and $\phi_2(p)>\mu_2$.
(2) for every $p\in S2(f)$ we have: $\phi_1(p)>\mu_1$
and $\phi_2(p)<\mu_2$.

For $\delta>0$ denote by $D1_\delta (v)$, resp. by
$D1_\delta(-v)$,  the intersection
with $V_0$, resp. with $V_1$,
of $\cup_{p\in S1(f)} D_\delta (p,v)$,
resp. of $\cup_{p\in S1(f)} D_\delta (p,-v)$.
By abuse of notation
the intersection of
$\cup_{p\in S1(f)} D(p,v)$
with $V_0$ will be denoted by
$D1_0(v)$.
Denote $D2_\delta(-v)\cup\stexp {(-v)} (B_a)$ by $\Delta(\delta,-v)$
and
$D1_\delta(v)\cup\st v (A_b)$ by $\nabla(\delta,v)$.
The similar notation like
$D2_\delta(-v),~ \Delta(0,-v),
\mathbreak
    D2_0(v)  $
etc.
are now clear without special definition.
For
 $\delta>0$    sufficiently small  we have
$$\gathered
\forall p\in S1(f): D_\delta(p)\subset\phi_1^{-1}(]a,\mu_1[)
\aand D_\delta(p)\subset
\phi_2^{-1}(]\mu_2, b[)\\
\forall p\in S2(f): D_\delta(p)\subset\phi_1^{-1}(]\mu_2, b[)
\aand D_\delta(p)\subset
\phi_2^{-1}(]a, \mu_2[)
\endgathered\tag\text{D1}
$$
It is easy to prove that for $\delta>0$
sufficiently small we have:
$$\nabla(\delta,v)\subset \Int A_a,\quad
\Delta(\delta,-v)\subset\Int B_b,\quad
\Delta(\delta,-v)\cap
D1_\delta(-v)=\emptyset
\tag\text{D2}
$$
Fix some $\delta>0$ satisfying (D1) and (D2).
\vskip0.1in\noindent
\it           Homomorphism \rm
$\h H(v;\mu',\mu; U): H_*(\qu (V_1\setminus B_b),\qu(A_b))\to
H_*(\qu(V_0\sm B_a),\qu(A_a))$
\vskip0.1in
Let $0\leqslant\mu '<\mu\leqslant\delta$.
Let $U$ be any subset of $V_1$ such that
$$\Delta (0,-v)\subset U
\subset B_b
\aand U\cap D1_{\delta}(-v)=\emptyset
$$
(for example $U=\Delta (\delta ,-v)$ will do).
 Denote by
$\h{\HH} (v;\mu',\mu ; U)$
the following sequence of homomorphisms
$$\gather
H_*(\qu(V_1\setminus B_b),\qu(A_b))
@>        \h I_*   >>
H_*\big( \qu(V_1\setminus U)~,~
\qu(A_b\cup D1_\mu (-v))\big)
@> {\text{Exc} ^{-1}} >>\\
H_*\big( \qu(V_1\setminus (U
\cup D1_{\mu'}(-v)))~,~
\qu\big((A_b\cup D1_\mu (-v))\setminus D1_{\mu'}(-v)\big)
@> \stexp {\h v}_* >>            \\
  H_*(\qu(V_0\setminus B_a),\qu(A_a)\big).
\endgather
$$
Here $\h I$
is the corresponding inclusion.  Note that the
 last arrow is well defined
since
$\stexp {(-v)} (B_a)\subset U\aand D_0(-v)\cap V_1\subset
D2_0(-v)\cup D1_{\mu'}(-v).$
All the three arrows are homomorphisms
of right $\ZZZ H$-modules;
first two - by the obvious reasons, the last
because $\stexp {\h v}$ commutes with
the right action of $H$.

The composition
$\stexp {\h v}_*\circ {\text{Exc}}^{-1} \circ \h I_*$
of this sequence will be denoted by
$\h H(v;\mu',\mu;\UU)$. The reasoning similar to
[9], page 25
 shows that this
 homomorphism does not depend on
 the choice
of $U$,neither
on the choice of $\mu',\mu$ or
$\delta$, or on the choice
of presentation $S(f)=S1(f)\sqcup S2(f)$
(if there is more then one such presentation).
Therefore this homomorphism is well determined by
$v$, the ranging system $\ran$
and the covering $\QQ:\h W\to W$.
We shall denote it by $\h H (v)$.
The proof of properties 1) and 2) of $\h H(v)$
is similar of the proof of [9, 4.7 (1,2)] and will be omitted.
\quad$\square$

\subhead{E. Equivariant ranging systems and the proof
of Theorem 1.3}
\endsubhead
We return here to the terminology of Subsection C.1 of
the Introduction.
We begin by some algebraic preliminaries. Let $M,N$ be
right $\ZZZ H$-modules and $f:M\to N$ be a
homomorphism of abelian groups. We say that $f$ is
$\theta$-\it semilinear \rm , if we have
$f(xh)=f(x)\t h\t^{-1}$ for every $x\in M$. If $M,N$
are free finitely generated $\ZZZ H$-modules with bases
$(e_j),(d_i)$, one can associate to each
$\t$-semilinear homomorphism $f:M\to N$ a matrix
$M(f)=(m_{ij})$ by the following rule:
$f(e_j)=\sum\limits_i d_im_{ij}$.

If $M$ is a right $\ZZZ H$-module and $f:M\to\ZZZ$ is a
\ho of abelian groups, then we shall say that
$f$ is of finite type, if for every $m\in M$
the set of $h\in H$, such that
$f(mh)\not=0$ is finite. If $f:M\to\ZZZ$ is a \ho of
a finite type, then we define a \ho $\tilde f:M\to\ZZZ H$
of $\ZZZ H$-modules by
$\tilde f (m)=\sum\limits_{h\in H} f(mh)h^{-1}$.

Returning to Morse maps,  we assume moreover that
$f:M\to S^1$  belongs to an indivisible cohomology class
in $H^1(M,\ZZZ)$.
Further, denote by $\pi:\bar M\to M$
the (unique) infinite cyclic covering, such that
$f\circ \pi\sim 0$.
The universal covering $p:\wi M\to M$ factors as
$p=\pi\circ\QQ$ where $\QQ:\wi M\to\bar M$
is a covering with structure group $H=\Ker\xi$.
 Let $u$ be an $f$-gradient
and let $\rans$ be a $t$-equivariant
ranging system for $(F,u)$ (see [9,Def. 4.14] for definition).

For $\nu,\mu\in\Sigma,\nu<\mu$ denote by
$\h H_{[\mu,\nu]}(u)$ the homomorphism
$\h H(u\mid F^{-1}([\nu,\mu]))$, associated
by virtue of
Proposition 1.5 to the ranging system
${\{(A_\sigma , B_\sigma)\}_
{\sigma\in\Sigma,~\nu\leq\sigma\leq\mu}}$
and the covering $\QQ$ (restricted to the cobordism
$F^{-1}([\nu,\mu])$). Denote by
$\h H_{[\mu,\mu]}(u)$ the identity homomorphism
of
$H_*\big(\qu(F^{-1}(\mu)\setminus B_\mu),\qu(A_\mu)\big)$
to itself.

It follows from the construction that for every
$g\in G$ with $\xi(g)=k\in\ZZZ$
we have
$\h H_{[\mu+k,\nu+k]}(u)=
R(g)\circ\h H_{[\mu,\nu]}(u)\circ R(g^{-1})$.
We have also $\h H_{[\nu,\theta]}(u)\circ\h H_{[\mu,\nu]}(u)=
\h H_{[\mu,\theta]}(u)$.
For $\nu\in\Sigma$ denote $R(\theta^{-1})\circ
\h H_{[\nu,\nu-1]}(u)$ by $\h h_\nu(u)$. It is a
$\theta$-semilinear endomorphism
of
$H_*\big(\qu(F^{-1}(\nu)\setminus B_\nu),\qu(A_\nu)\big)$.
We have obviously
$\h H_{[\nu,\nu-k]}(u)=R(\theta^k)\circ (\h h_\nu (u))^k$.
The next lemma
follows directly from 1.5.

\proclaim{Lemma 1.6}
Let $\mu,\nu\in\Sigma,\nu\leq\mu$; let $k\in\NNN$.
Let $N$ be oriented-lifted submanifold of
$F^{-1}(\mu)\setminus B_\mu$ such that $N\setminus\Int A_\mu$
is compact. Let $L$ be a cooriented-lifted compact
submanifold of $F^{-1}(\nu)\setminus A_\nu$. Assume that
$\dim N+\dim L=\dim M-1$.
Then:
\roster\item
$N'_k=
 \stind u\mu{\nu -k} (N)$ is an oriented-lifted
submanifold of $\fpr F{\nu -k}\setminus B_{\nu -k}$
 such that
\break
$ N'_k\setminus\text{\rm{ Int}}~ A_{\nu-k}$
is compact.
If $N'_k\pitchfork L t^k$, then
$N'_k\cap L t^k$ is finite and
$\h N'_k ~\sharp ~ \h L\theta^k=
 \h i^*( ]\h L[)\big((\h h_{\nu}(u))^k([\h N_0'])\big)$,
 where $\h i$ stands for the inclusion map
\newline
$\big(\qu(F^{-1}(\nu)\setminus B_\nu) , \qu(A_\nu )\big)
\hookrightarrow
\qu\big(F^{-1}(\nu)  ,\qu(F^{-1}(\nu)\setminus \h L)\big)$.

\item For every $f$-gradient $w$, sufficiently close
 to $u$ in $C^0$-topology, $\rans$ is also
 a $t$-equivariant ranging system for $(F,w)$ and
$\h h_{\nu}(u)=\h h_{\nu}(w),
\break
\h H_{[\mu,\nu]}(u)= \h H_{[\mu,\nu]}(w). \qquad\square$
\endroster
\endproclaim

\it
Proof of Theorem 1.3
\rm\quad
It is easy to see that it suffices to prove our theorem for
the case of indivisible homotopy class
$[f]\in H^1(M,\ZZZ)$ and we make this assumption up to the end
of this subsection.

Fix first two points $x,y\in S(f), \ind x=\ind y+1$.
Recall that we have fixed a lifting $\wi x\in\wi M$
for every $x\in S(f)$; denote $\QQ(\wi x)$ by $\bar x$.
We can assume that
$F(\bar y)<F(\bar x)\leq F(\bar y)+1$.
Denote $\dim M$ by $n$; denote $\ind x$ by $l+1$,
then $\ind y=l$.
Choose some set $\Sigma$ of regular values
of $F$, satisfying $({\Cal S})$ of [9, Def. 4.14].

Denote by $\theta$ the maximal element of $\Sigma$ with
$\theta< F(\bar x)$ and by $N(v)$ the
 intersection $D(\bar x,v)\cap F^{-1}(\theta)$;
$N(v)$ is an oriented submanifold of $F^{-1}(\theta)$,
diffeomorphic to $S^l$.
Denote by $\eta$ the minimal element of $\Sigma$,
satisfying $\eta>F(\bar y)$;
then $\eta\leq\theta<\eta+1$.
Denote by $L(-v)$ the
 intersection $D(\bar y,-v)\cap F^{-1}(\eta)$;
$L(-v)$ is a cooriented submanifold of $F^{-1}(\eta)$,
diffeomorphic to $S^{n-1-l}$.
Denote by $W$ the cobordism $F^{-1}([\eta,\eta+1])$.
Note that $\bar x\in \Wkr$. Denote
$F^{-1}(\eta)$ by $V_0$,
$F^{-1}(\eta+1)$ by $V_1$,
$\Sigma\cap[\eta,\eta+1]$ by $\Lambda$.

Denote by $\GG t_1(f;x,y)$ the subset of $\GG t(f)$,
consisting of all the $f$-gradients $v$, such that there is
an equivariant ranging system
$\rans$ for $(F,v)$ satisfying
$$\gather
N(v)\cap B_\theta =\emptyset,
\quad L(-v)\cap A_\eta=\emptyset, \tag1.1\\
\left\{
\gathered
(F^{-1}(\eta)\setminus B_\eta,A_\eta)
\text{\rm~~ has a homotopy type
of a finite CW-pair} \\
\text{\rm      having only
l-dimensional cells.~~}
\endgathered  \right. \tag1.2
\\
\endgather
$$
Now we shall prove 3 properties of the set
$\GG t_1(f;x,y)$.
\subsubhead{(1).
 $\GG t_1(f;x,y)$ is an open and dense
subset of $\GG t(f)$ with respect to
 $C^0$ topology}
\endsubsubhead

This is proved exactly in the same way
as the open-and-dense property
of
$\GG t_0(f;x,y)$ in [9,p.29,30].

\subsubhead{(2).
 If $v\in\GG t_1(f;x,y)$ then $\wi n(\wi x,\wi y;v)$ satisfies
the condition $(\LL)$}
\endsubsubhead

The liftings $x\mapsto\wi x$ and $y\mapsto\wi y$ define a lifting
$\h N(v)$ of $N(v)$ to $\qu\big(F^{-1}(\theta)\big)$
and 
$\h L(-v)$ of $L(-v)$ to $\qu\big(F^{-1}(\eta)\big)$.
 The $\ZZZ H$-module
$\HH=H_l((\qu(F^{-1}(\eta)\setminus B_\eta) , \qu(A_\eta ))$
is free; choose some basis $e_1,...,e_m$
of this module.
The homomorphism
$\h h_\eta(v)$ of this module is $\theta$-semilinear.
Denote by ${\Cal B}=(b_{ij})$ its matrix, and
denote by $\Cal A$ the matrix $(b_{ij}\t)$. Let
$a_{ij}^{(s)}$ be the coefficients of $\AA^s$.
Consider the element $\xi=\h H_{[\t,\eta]}(v)([\widehat N (v)])$
of $\HH$; let $\xi=\sum e_iX_i$ with $X_i\in\ZZZ H$.
Consider $\beta=\h i^*(]\h L(-v)[)$ as a
 homomorphism of $\HH$ to $\ZZZ$. It is of finite type;
denote $\tilde\beta ( e_j)$
by $Y_j$; then $Y_j\in\ZZZ H$.
We claim that
$$
\wi n(\wi x,\wi y;v)=\sum
\Sb
s\geq 0\\
1\leq i\leq m\\
1\leq j\leq m
\endSb
Y_i\cdot a^{(s)}_{ij}\cdot X_j
\tag*
$$
To prove it write
$\wi n(\wi x,\wi y;v)=\dsize\sum\limits_{s\geq 0}
\bigg(\dsize\sum\limits_{h\in H}
\nu(\wi x,\wi yh\theta^s )
\cdot h\theta^s\bigg)$
(here $\nu(\wi x,\wi yh\theta^s )$ stands for the algebraic
number of $(-v)$-trajectories, joining
$\tilde x$ with $\tilde yh\theta^s$;
note that since $F(\bar x)\leq F(\bar y)+1$,
there are no $(-v)$-trajectories joining
$\wi x$ to $\wi yg$ if $\xi (g)>0$).

To make the following computation more easy to comprehend,
we make the following terminology conventions (valid only here).
The homomorphism $\h h_\eta(v):\Cal H\to \Cal H$
will be denoted by $\mu$. We  identify the cohomology classes
in 
$H^*(\qu\big(F^{-1}(\nu)  ,\qu(F^{-1}(\nu)\setminus \h L)\big)$
with their images in
$H^*\big(\qu(F^{-1}(\nu)\setminus B_\nu) , \qu(A_\nu )\big)$
(thus suppressing $\h i^*$ in the notation).
We have:
$$\wi n(\wi x,\wi y;v)=
  \sum\limits_{s\geq 0}\bigg(   \sum\limits_{h\in H}
\big(
]\h L\cdot h[
\big( \mu^s(\xi))\big)h \bigg)\theta^s$$
(by 1.6).
 The latter expression equals to 
$$\sum\limits_{s\geq 0}\bigg(   \sum\limits_{h\in H}
\big(
]\h L[
( \mu^s(\xi)\cdot h^{-1})\big)h \bigg)\theta^s
=
\sum\limits_{s\geq 0}\bigg(   \sum\limits_{h\in H}
\b(\m^s(\xi)\cdot h^{-1}\big)h
 \bigg)\theta^s
=
\sum\limits_{s\geq 0}\widetilde\b(\m^s(\xi))\cdot\theta^s$$

To obtain from this expression the formula (*) we need only
a lemma, allowing to calculate $\m^s(\xi)$ in terms of the coordinates
of $\xi$ and the matrix of $\m$ (the expression is slightly
different from the standard one in linear algebra since $\m$
is $\t$-\it semilinear \rm ).

\proclaim{Lemma 1.7}
Let $F$ be a free $\ZZZ H$-module with a basis
$e_1,...,e_m$ and $\mu:F\to F$ be a $\theta$-semilinear
homomorphism of $F$.
Let $m_{ij}$  be its matrix. Denote by $M$ the $m\times m$
matrix $(m_{ij}\cdot\theta)$.
Let $\xi\in F,\xi=\sum e_j\xi_j$.

 Then for every natural
$s\geq 0$ we have
$$\mu^s(\xi)=\sum_{i,j}e_i[M^s]_{ij}\xi_j\theta^{-s}$$
\endproclaim
\demo{Proof}
Induction in $s$. We have
$$\gather
\mu^{s+1}(\xi)=\mu(\mu^s(\xi))=
\mu\bigg(\sum_ie_i\bigg(\sum_j[M^s]_{ij}\xi_j\theta^{-s}\bigg)\bigg)=
\sum_i\mu(e_i)\cdot\bigg(\theta\sum_j[M^s]_{ij}\xi_j\theta^{-s-1}\bigg) \\
=\sum_i\bigg(\sum_ke_km_{ki}\bigg)\cdot
\bigg(\theta\sum_j[M^s]_{ij}\xi_j\theta^{-s-1}\bigg)=
\sum_ke_k\bigg(\sum_i(m_{ki}\theta)\cdot[M^s]_{ij}\bigg)
\xi_j\theta^{-s-1}.
\quad\square
\endgather
$$                                                  \enddemo

Now substitute the expression for $\m^s(\xi)$ into the above
formula, and get (*).

\subsubhead{(3).
 Let $v\in\GG t_1(f;x,y)$.
Let $U$ be a neighborhood of $S(f)$. Then there is $\epsilon>0$
such that for every $w\in\GG t_0(f;x,y)$
with
$\Vert w-v \Vert<\epsilon$
and $w\vert U=v\vert U$
we have:
$\wi n(\wi x,\wi y;v)=\wi n(\wi x,\wi y;w)$}
\endsubsubhead

Let $w$ be an $f$-gradient, sufficiently close to
$v$. Then $\rans$ is still a $t$-equivariant
ranging system for $(F,w)$, satisfying (1.1), (1.2).
It is not difficult to see 
that $[\h N(w)]=[\h N(v)], [\h L(-w)]=[\h L(-v)]$.
Then Lemma 1.5 together with the formula (*)
above finishes the proof.
\vskip0.2in
Set now $\GG t_1(f)$ to be the intersection of
all the $\GG t_1(f;x,y)$ where
$x,y\in S(f),\ind x=\ind y+1$, and the
Theorem 1.3 is proved.
\quad$\square$

\remark{Remark 1.8}
There is an obvious analog of Theorem 1.3 for any regular covering
$\phi:\h M\to M$ such that $\phi^*[\o]=0$.
\endremark

\newpage

\head
\S 2. Morse forms with arbitrary cohomology classes
\endhead

\subhead
A. Algebraic preliminaries
\endsubhead

We shall need some lemmas about the ring
$\gsq {\ZZZ^m}$ and its completions
and localizations.

Let $\eta:\ZZZ^m\to\RRR$ be a non-zero \ho . We
extend it to a linear map $\RRR^m\to\RRR$, which will be
denoted by the same letter. We say that a set $Z\subset\RRR^m$
is an \it $\eta$-cone   \rm , if there is a compact convex nonempty set
$K\subset\eta^{-1}(-1)$ such that
$Z=\{\l z\mid\l\in\RRR,\l\geq0,z\in K\}.$
We say, that $Z$ is $(\xi,\eta)$-cone if $Z$ is $\xi$-cone
and $\eta$-cone. We say that a set $Z\subset\RRR^m$ is an
\ic {$\eta$} if there are $e_1,...,e_k\in\ZZZ^m$, such that
\roster\item
$\rk (e_1,...,e_k)=m$
\item $\eta(e_i)<0$
\item $Z=\{\l_1e_1+...+\l_ke_k\mid\l_i\geq0\}$
\endroster
We shall also write $Z=Z\langle e_1,...,e_k\rangle$.
We say, that $Z$ is an \ic {$(\xi,\eta)$} if the vectors
$e_1,...,e_k$ above satisfy $\xi(e_i)<0,\eta(e_i)<0$ for all $i$.
Note that an \ic {$(\xi,\eta)$} is a $(\xi,\eta)$ - cone.

\proclaim{Lemma 2.1}
Let $Z$ be an $(\eta_1,\eta_2)$-cone. Then there is an integral
$(\eta_1,\eta_2)$-cone $Z_0\supset Z$.
\endproclaim
\demo{Proof}
We assume that $\eta_1, \eta_2$ are linearly independent;
the other case is considered similarly.
Denote $\eta_1^{-1}(-1)$ by $H$, and the set
$H\cap\{\eta_2(x)<0\}$ by $H_0$;
then $H_0$ is an open halfspace of $H$, containing $Z\cap H$.
Denote by $\LL$ the set $\{\l x\vert \l\geq 0,x\in\ZZZ^m\}$.
$\LL$ is everywhere dense in $H$ and in $H_0$. It is
not difficult to prove that there is a finite subset
$\LL_0\subset\LL$ such that $\LL_0\subset H_0$ and
$\langle \LL_0\rangle\supset Z\cap H$. We can choose $\LL_0$ so that
$\rk \LL_0=m$ and the lemma is proved.
\quad$\square$\enddemo

\proclaim{Lemma 2.2}
Let $\xi:\RRR^m\to\RRR$ be a non-zero linear form.
Let $\e>0$. Then there is a finite set $I$ of linear forms
$\t_i:\RRR^m\to\RRR, i\in I$ with
$\Vert \t_i\Vert<\e$ such that
\roster\item
the set $\G=\{x\in\RRR^m\mid (\xi+\t_i)(x)\leq 0\}$
is a $\xi$-cone.
\item
There is an \ic {$\xi$} $\G_0$ such that for every family
$\{A_i\}_{i\in I}$ of real numbers there is
$b\in \ZZZ^m$ with the property:
$$\{x\in\RRR^m\mid (\xi+\t_i)(x)\leq A_i\}\subset \G_0+b\tag*$$
\endroster
\endproclaim
\demo{Proof} 1)
Pick any $m$ linearly independant linear forms
$\a_1,...,\a_m:\RRR^m\to\RRR$
with $\Vert\a_i\Vert<\min(\e,\Vert\xi\Vert/2)$.
I claim that the finite family
$\{\a_1,-\a_1,...,\a_m,-\a_m\}$
of linear forms satisfy our conclusions.
Note first that
$$Z=\xi^{-1}(-1)\cap\{x\mid
\forall i~:~
(\xi+\a_i)(x)\leq 0,(\xi-\a_i)(x)\leq 0\}$$
is non empty and compact.
(Indeed, let $x$ be a vector such that
 $\vert x\vert =1$ and $\xi(x)=\Vert\xi\Vert$.
Then $a=-\frac x{\xi(x)}\in Z$. Further, if $Z$ is not bounded, then
there is a sequence $x_n\in Z$ such that
$\vert x_n\vert\to\infty$. Consider the sequence
$x_n/\vert x_n\vert$. We can assume that it converges to some
 $v$,$\Vert v\Vert=1$. Since $\xi(x_n)=-1$,
we have $\xi(v)=0$. Further, for every $i$ we have
$(\xi+\a_i)(v)\leq0,(\xi-\a_i)(v)\leq0$, \th
$\a_i(v)=0$ and $v=0$.)
Further, $\G\subset\xi^{-1}(]-\infty,0])$ and
the intersection $\G\cap\xi^{-1}(0)$ consists of $0$.
\Th $x\in\G,x\not= 0$ implies $x=\l y,y\in Z$.
2)
Choose a vector $x_0\in\ZZZ^m$ such that $\xi(x_0)<0$ and for every
$1\leq i\leq m$ we have
$(\xi+\a_i)(x_0)<0,(\xi-\a_i)(x_0)<0$
(such a vector exists; it suffices to note that
$(\xi\pm\a_i)(a)=-1\mp\frac {\a_i(x)}{\xi(x)}<0$, and to
approximate $a$ by an element $y_0\in\QQQ^m$.)

Denote $\vert (\xi\pm\a_i)(x_0)\vert$ by $\b^{\pm}_i>0$.
Now $(\xi\pm\a_i)(x)\leq A_i$ implies that
$(\xi\pm\a_i)(x+px_0)\leq A_i-p\b^{\pm}_i.$
\Th if $p$ is sufficiently big, we obtain (*) with
$b=-px_0$, or, equivalently,
$\{x\in\RRR^m \mid (\xi+\t_i)(x)\leq A_i\}\subset\G+b$.
Then apply 2.1 to obtain an integral $\xi$-cone $\G_0$, such that
$\G\subset\G_0$.
\quad$\square$
\enddemo

Next we shall recall some facts from [6].
Let $\eta:\ZZZ^m\to\RRR$ be a \ho
 and $Z=Z\langle e_1,...,e_k\rangle$ be an integral
$\xi$-cone. Denote $\ZZZ^m\cap Z$
by $[Z]$, and the set
$\{n_1e_1+...+n_ke_k \mid n_i\in N\}$ by $[\wi Z]$.
Note that $[Z]$ and $[\wi Z]$ are submonoids of
$\ZZZ^m$, $[\wi Z]\subset [Z]$, and $[Z]$ is finitely generated
over $[\wi Z]$. Consider the group rings
$\gsq Z$ and $\gsq {\wi Z}$ of the monoids $Z$ and $\wi Z$;
since $\gsq {\wi Z}$ is an image of the ring
$\ZZZ[t_1,...,t_k]$, it is noetherian. The ring
$\ZZZ[Z]$ is finitely generated as a left module over
$\ZZZ[\wi Z]$, \th it is also noetherian.

Consider the submonoid $\{x\in Z\mid x\not=0\}$ of $Z$
and denote its group ring by $\frak m_Z$.
It is an ideal of $\zd$. Denote by $\zd\sphat$
the $\frak m_Z$-completion of $\zd$.
Denote by $S_Z$ the multiplicative set $1+\frak m_Z\subset\zd$.
Then $(S^{-1}_Z\zd)\sphat=(\zd)\sphat$
(see [2], Ch 3,\S 3, Prop. 12).
The ring $(\zd)\sphat$ is easily identified with
the ring of all the elements $\l\in(\zd)\sphat~\sphat$ such that
$\supp\l\subset Z$ (this
latter ring is obviously a subring of $\zd_\eta^-$).
The ring $(\zd)\sphat$ is  faithfully flat
over $S^{-1}_Z\zd$ (ibid. Prop. 9).
\Th , if $P,Q\in\zd$ and there is $x\in(\zd)\sphat$
such that $P=Qx$, then $x\in S^{-1}_Z\zd$ (ibid. Ch 1,\S 3,p.7).

Set $\s_Z=\{t^I\mid I\in Z\}$; it is a multiplicative subset of $\zd$.
It is easy to see that $\s_Z^{-1}\zd=\zzz$; and
$\s_Z^{-1}(\zd\sphat~)=
\{S\in(\zd)\sphat~\sphat\mid \exists x\in\ZZZ^m:\supp S\subset Z+x\}$.
The faithful flatness property cited above implies immediately that
if $P,Q\in\zzz$ and there is $x\in\s_Z^{-1}(\zd\sphat~)$
such that $P=Qx$, then $x\in\s^{-1}_ZS_Z^{-1}\zd$.

For a linear form $\eta$ denote by $S_\eta$ the multiplicative subset
$\{1+Q\mid \supp Q\subset
\mathbreak
\eta^{-1}(]-\infty,0[)\}$.

\proclaim{Lemma 2.3}
Let $\a:\ZZZ^m\to\ZZZ$ be an indivisible \ho .
Denote $\Ker\eta$ by $H$.
Let $A$ be $(k\times k)$-matrix, such that $a_{ij}\in\zzz_{(-1)}$.
Let $\xi=(\xi_1,...,\xi_k),\eta=(\eta_1,...,\eta_k)$ be vectors in
$(\ZZZ H)^k$. Denote by $a_{ij}^{(s)}$ the $(ij)$-coefficient
of the matrix $A^s$.

Then $(1-det A)\sum_{s\geq0}(\sum_{i,j}\xi_ia_{ij}^{(s)}\eta_j)$
belongs to $\zzz$.
\endproclaim
\demo{Proof}
It suffices to prove that every coefficient of the $(k\times k)$-
matrix
$(1- det A)\bigg(\sum\limits_{s\geq 0}A^s\bigg)$
belongs to $\zzz$.
Consider the matrix $1-A$. It is invertible in the ring
$S_\a^{-1}\zzz$ and the Cramer rules imply
$(1-detA)(1-A)^{-1}\in \text{\rm Mat}(\zzz)$.
On the other hand $(1-A)^{-1}=    \sum\limits_{s\geq 0} A^s$.
\quad$\square$\enddemo

\subhead
B. Preliminaries on Morse forms and their gradients
\endsubhead

In this subsection we assume the terminology of \S 1 of [9].
Moreover, we assume the definition of $\o$-chart-system
and, respectively, of $\o$-gradient (where $\o$ is a Morse
form), which are completely similar to those of $f$-chart-
system and $f$-gradient, see Definition 1.1 of [9].

So let $\omega$ be a Morse form on a closed connected
manifold $M$.
Let $p:\widetilde M\to M$ be the universal covering
and let $\widetilde F:\widetilde M\to\RRR$
be a Morse function such that $d\wi F=p^*\o$.
Note that $\wi F(xg)=\wi F(x)+\acc \o (g)$,
where $g\in\pi_1M$.
Choose a riemannian metric on $M$. Then $\wi M$
obtains a $\pi_1M$-invariant riemannian metric.
If $v$ is an $\omega$-gradient
we denote by the same letter $v$ its liftings to
$\widetilde M$ and $\bbar M$, since
\noconf .
The length of a curve $\gamma$ will be denoted by
$l(\gamma)$.
The following simple and useful lemma is known since the early
80s. I knew it from J.-Cl.Sikorav.

Assume that for every $x\in S(f)$ a lifting $\widetilde x$ of $x$ to
$\widetilde M$ is fixed.

\proclaim{Lemma 2.4} Let $v$ be an $\omega$-gradient. Then there
are constants $A,B>0$ such that:

For every $g\in\pi_1M$ and every $(-v)$-trajectory $\gamma$,
joining $\widetilde x$ and $\widetilde y\cdot g$
we have:
$l(\gamma)\leq A-B\acc {\omega}(g)$.
\endproclaim
\demo{Proof}
Choose any $\delta>0$ less than the injectivity radius of $M$
and less than $\min\limits_{p,q\in S(\o)}\rho(p,q)$.
Then any nontrivial piecewise smooth loop in $M$
is longer than $\delta$. \Th for any
$x\in\widetilde M, 1\not=g\in\pi_1M$, any piecewise
smooth path in $\widetilde M$, joining a point
of $D(x,\delta/3)$ with a point of $D(xg,\delta/3)$
has the length $\geq\delta/3$.
Also for $x, y \in\wi M$ with
$p(x),p(y)\in S(\o)$ any piecewise smooth path
joining joining a point
of $D(x,\delta/3)$ with a point of
$D(y,\delta/3)$ has the length $\geq\d/3$.

Now let
$\{\Phi_p:U_p\to B^n(0,r)\}_{p\in S(\o)}$
be an $\omega$-chart-system such that
$\GG(\overline U_p,\Phi_p)\leq C$ for every $p\in S(\o)$
and that $rC<\d/12$.\footnote{Recall from [9,p.3] that
for a chart $\Phi:U\to V\subset\RRR^n$
of $M$ and $x\in U$
we denote by $\GG(x,\Phi)$
the number
$\sup\limits_{x\in M, h\in T_xM, h\not=0}
(\max \big(\vert h\vert_\rho/\vert\Phi_*h\vert_e,
\vert\Phi_*h\vert_e/\vert h\vert_\rho)\big)$,
where $\rho$ stands for the metric on $M$ and $e$ for
the euclidean metric in $\RRR^n$.} This condition implies
in particular that $U_p\subset D(p,\d/12)$.

Choose some liftings $\widetilde U_p$ of neighborhoods
$U_p$, extending $x\mapsto\widetilde x$.
Let $D>0$ be less than $\min\limits_{x\notin \cup U_p}
\omega(v(x))$. Denote $\Vert v\Vert$ by $E$.

Now let $\gamma$ be a $v$-\tr , joining $\wi x$ with $\wi yg$,
and let
$A_0=\wi x, A_1,...,A_n,\wi y=A_{n+1}$ be the points in
$p^{-1}(S(\omega))$ such that $\gamma$ intersects
$\wi U_{A_i}$. Then $(n+1)\cdot\delta/3\leq l(\gamma)$.
The length
of the part of $\gamma$ inside of
$\bigcup\limits_{i=0}^{n+1}\wi U_{A_i}$ is not more than
$2rC(n+1)$.

Now let $t_i$, resp. $\tau_i$ be the moment when $\gamma$ enters
$\wi U_{A_i}$, resp. quits $\wi U_{A_i}$.
We have
$$\int_{\tau_i}^{t_{i+1}}\omega(\overset\cdot\to\gamma)dt=
 -\int_{\tau_i}^{t_{i+1}}\omega(v)dt=
\wi F(\g(t_{i+1}))-\wi F(\g(\tau_i)).$$

Therefore the total time which $\gamma$
can spend outside
$\bigcup\limits_p\wi U_{A_i}$ is not more
than
\break
$\vert \wi F(\wi yg)-\wi F(\wi x)\vert/D$, and the length
of the corresponding part of the curve is $\leq
E/D\vert \wi F(\wi yg)-\wi F(\wi x)\vert.$
Since $\gamma$ joins $\wi x$ to $\wi y g$,
the last expression
 is $\leq E/D\cdot (\wi F(\wi x)-\wi F(\wi y)-\acc \omega(g))$.
\Th
$$l(\gamma)\leq 2rC(n+1)+E/D (\wi F(\wi x)-\wi F(\wi y)) -
E/D\acc \omega (g) \leq
\frac {6rC}{\d}l(\g)  +\frac A2 -\frac B2\acc \omega(g) $$
where $A$ is chosen such that
 $A/2\geq \vert \wi F(\wi x) -\wi F(\wi y)\vert$
for every $x,y\in S(f)$ and $B=2E/D$. Then the inequality
$l(\g)\leq A-B\acc \omega(g)$ easily follows.
 \quad$\square$
\enddemo

Let $\o$ be a Morse form, let $\alas$
be an $\omega$-chart-system. Choose a basis $a_1,...,a_m$ in
$H_1(M,\ZZZ)/\text{Tors}$. Choose and fix closed 1-forms
$\lambda_1,...,\lambda_m$ on $M$, such that
$\langle[\lambda_i],a_j\rangle=\delta_{ij}$
and $\supp\lambda_i \cap\overline U_p=\emptyset$
for every $i$ and every $p\in S(\omega)$.
(To prove that we can satisfy the second condition, let $\theta$
be any closed 1-form. Let
$\{\wi\Phi_p: U'_p\to B^n(0,r'_p)\}_{p\in S(\omega)}$
be some standard extension of
$\alas$ and let
$\phi_p$
be a $\smo$ function which equals to 1
in a neghborhood of $\overline U_p$ and ~~
$\supp\phi_p\subset U'_p$.
Let $F_p$ be a function on $ U'_p$, such that
$dF_p=\theta$. Consider the form $\theta'=
\theta-\sum\limits_{p\in S(\omega)}d(\phi_pF_p)$. We have
$\sq {\theta'}=\sq \theta$ and $\theta(x)=0$
in every $\overline U_p$.)

For $\vec\e=(\e_1,...,\e_m)\in\RRR^m$ we denote by
$\vec\e\cdot\vec\lambda$ the form
$\sum\limits_{i=1}^m \e_i\lambda_i$,
by $\omega_{\vec\e}$ the form $\omega+\vec\e\cdot\vec\lambda$.
For $\e>0$ set
$$\Omega_\e =\{\omega_{\vec\e}\mid \vert\vec\e~\vert
=\max_i\vert\e_i\vert\leq\e\}.$$

We shall say that
$\O_\e$ is a \it Morse family \rm, if for every
$\vec \e$ with $\vert \vec\e\vert\leq\e$ the form
$\o_{\vec\e}$ is a Morse form and $S(\o_{\vec\e})=S(\o)$.
Let $\O_\e$ be a Morse family, and $v$ be a vector field. We say, that
$v$ is an $\O_\e$-\it gradient \rm, if $v$ is an $\o_\e$-gradient
for each $\o_{\vec\e}\in\O_\e.$

\proclaim{Lemma 2.5}
\roster\runinitem
There is $\e>0$, such that $\O_\e$ is a Morse family.
\item Let $v$ be an $\o$-gradient. Then there is $\e>0$ such that
$v$ is an $\O_\e$-gradient
\item Let $v$ be an $\O_\e$-gradient. Then there is $\d>0$
such that every $\o$-gradient $u$ with $\Vert u-v\Vert<\delta$
is an $\O_\e$-gradient.
\endroster
\endproclaim
\demo{Proof}
1) Denote $\sup\Sb 1\leq i\leq m\\ x\in M \endSb      \Vert \l_i(x)\Vert$
by $\l$ and $\min\limits_{ x\in M\setminus\cup_pU_p}\Vert \o(x)\Vert $
by $\eta$. Then $\e=\frac {\eta}{2m\l}$ will do.

2) Denote $\min\limits_{ x\in M\setminus\cup_pU_p} \o(v)(x) $
 by $\bar \eta$ and
$\sup\Sb 1\leq i\leq m\\ x\in M \endSb  \vert \l_i(v)(x)\vert$
by $\bar \l$. Then $\e=\frac {\bar\eta}{2m\bar\l}$ will do.

3) Denote $\sup\Sb 1\leq i\leq m\\ x\in M \endSb     \Vert \l_i(x)\Vert$ by $\l$.
Denote by $Q$ the compact set
$(M\setminus\cup_pU_p)\times\sq {-\e,\e }^m$ and by
$F:Q\to\RRR$ the map
$F:(x,\nu_1,...,\nu_m)\mapsto\omega_{\vec\nu}(v)(x)$.
Since $\Im F\subset ]0,\infty[$, there is $\beta>0$ such that
$\Im F\subset [\b,\infty[$. Let $u$ be any $\o$-gradient such that
 $\Vert u-v\Vert\cdot\big(m\l\e
+\Vert\o\Vert\big)\leq\b/2$.

We claim that $u$ is an $\O_\e$-gradient. Indeed, note first
that $\o_{\vec\e}(u)(x)>0$ for any
$x\in (\cup_pU_p\setminus S(\o))$ and any
 $\vert\vec\e\vert\leq\e$.
Further, if $x\in (M\setminus\cup_pU_p)$, we have
$$
\vert\o_{\vec\e}(u)(x)-\o_{\vec\e}(v)(x)\vert=
\vert(\omega+\sum_i\e_i\l_i)(u-v)(x)\vert
\leq
\Vert u-v\Vert \cdot (\Vert\o\Vert+\e\l m)
\leq
\b/2
$$
therefore $\o_{\vec\e}(x)>0$. Finally, $u$ has a standard
form with respect to some $\o$-chart-system.

 The suitable restriction of this system
will be an $\o_{\vec\e}$-chart-system for any $\vec\e$.\quad
$\square$
\enddemo

Now we can define the incidence coefficients with respect to the universal
cover.
The preceding lemma implies that $v$ is an $\bar\o$-gradient
for some 1-form $\bar\o$ which cohomology class is rational.
\Th (see [7]) for every $g\in\pi_1(M)$
there is at most finite set
 of $(-v)$-trajectories joining $\wi x$ with $\wi yg$
if $\ind x=\ind y+1$. Choose orientations of descending discs.
For each such \tr we denote by $\e(\g)$ the sign of intersection of
$D(\wi x,v)$ with $D(\wi yg,-v)$ along $\gamma$. The element
$\sum\limits_{\g}\e(\g)$ is denoted by
$\nu(\wi x,\wi yg)$ and we set
$\wi n(\wi x,\wi y;v)=
\sum\limits_{g\in G}\nu(\wi x,\wi yg)g$. This
 is an element of the abelian
group $(\ZZZ[\pi_1M])\sphat~\sphat$
of all the formal linear combinations
(infinite in general) of the elements of $G$.

\proclaim{Lemma 2.6}
$\wi n(\wi x,\wi y;v)\in (\gsq {\pi_1M})^-_{\acc \o}$.
\endproclaim
\demo{Proof}
Recall that Novikov ring
$\ZZZ([\pi_1M])^-_{\acc \o}$
consists of all
$\l\in(\ZZZ[\pi_1M])\sphat~\sphat$
 such that
for every $c\in\RRR$ we have:
$\supp\l\cap{\acc \o}^{-1}([c,\infty[)$ is finite.
Then our lemma follows from Lemma 2.4.
\quad$\square$
\enddemo

\remark{Remark 2.7} Note that the analogs of Lemmas 2.4 
and 2.6
are obviously true for any regular covering $p':\widehat M\to M$
such that $(p')^*[\o]=0$.   \endremark

\subhead
C. Incidence coefficients with respect to
a free abelian cover
\endsubhead

Next we pass to free abelian coverings.
We assume here the terminology of Subsection C.2 of Introduction.
The universal covering factors then
 as $p=\PP\circ\QQ; \QQ:\widetilde M\to\bbar M$.
The epimorphism $\pi_1(M)\to H_1(M,\ZZZ)/
\text{~Tors}$ will be denoted by $Q$.
The deRham cohomology class $[\o]$ of $\o$ defines a
\ho
$H^1(M,\ZZZ)/\text{~Tors}
\to\RRR$, which will be denoted by the same letter
$[\o]$. Note that $\acc \o=\sq \o\circ Q$.
Since $\PP^*(H^1(M,\RRR))=0$, there is a Morse function
$\bbar F:\bbar M\to\RRR$, such that
$d\bbar F=\PP^*\o$, and we shall assume that
$\wi F=\bbar F\circ Q$.
 We have chosen  a riemannian metric on $M$; therefore
the manifold  $\bbar M$
obtain a riemannian metric, which is $\ZZZ^m$-invariant.
We have chosen a  basis $(a_1,...,a_m)$ in
$H_1(M,\ZZZ)/\text{\rm Tors}$. \Th this group is identified
 with $\ZZZ^m$, and the vector space
$H^1(M,\RRR)$ with the dual space of linear forms
$\RRR^m\to\RRR$. We choose the $L_1$-norm
in $\RRR^m$; then the dual space obtains the
$sup$-norm. (That is $\Vert\sum_i\a_ia_i\Vert=\sum_i\vert\a_i\vert$
and $\Vert\sum_i\b_ia_i^*\Vert=\max_j\vert\b_j\vert$,
where $\{a_i^*\}$ is the base dual to $\{a_i\}$.
)

\proclaim{Corollary 2.8}
Let $v$ be an $\o$-gradient.
There is such an $\e>0$, that every linear form $\eta:\RRR^m\to
\RRR$ with $\Vert [\o]-\eta\Vert\leq\e$ is a cohomology class
of a Morse form $\o(\eta)$ such that $v$ is an $\o$-gradient.
\endproclaim
\demo{Proof}
Every linear form $\eta:\RRR^m\to\RRR$ with $\Vert[\o]-\eta\Vert\leq\e$
can be written as $\eta=[\o]+\sum\eta_ia_i^*$ where
$\eta_i\in\RRR, \vert\eta_i\vert\leq\e$. Now let $\e>0$ be so small
that  $\O_\e$ is a Morse family and $v$ is an $\O_\e$-gradient
and set $\o(\eta)=\o+\sum\eta_i\l_i$.
$\square$\enddemo

For two critical points $x,y\in S(\om)$ and an $\om$-\gr
$v$ we set
$I(x,y;v)=
\mathbreak
\{g\in \pi_1(M)\mid $ there is a $(-v)$-\tr
$\g$, joining $\tilde x$ to $\tilde y\cdot g \}$.
If the set of $(-v)$-\trs joining $\tilde x$ to $\tilde y\cdot g$
is finite, we denote by $N(\tilde x,\tilde y,g;v)$ its cardinality.
(We identify here two \trs which differ by a parameter change.)

\remark{Remark 2.9} If $v$ is a good
$\om$-\gr and $x,y\in S(\om)$ with $\ind x=\ind y+1$, then
$\supp(\tilde n(\tilde x,\tilde y;v))\subset
I(x,y;v)$
and for every $g\in \pi_1M$ the set $N(\tilde x,\tilde y,g;v)$
is finite. \quad$\triangle$    \endremark

\proclaim{Lemma 2.10}
There is an \ic {$[\o]$} $\G$ and a vector $a\in\ZZZ^m$
such that
\break
 $Q(I(x,y;v))\subset\G+a$.
\endproclaim
\demo{Proof}
Note that if $v$ is a $\k$-gradient for some Morse form
$\k$, then
$Q(I(x,y;v))
\mathbreak
\subset \sq\k^{-1}(]-\infty, A])$
for some $A$. By Corollary 2.8 there is $\e>0$, such that
 every linear form $\eta:\RRR^m\to\RRR$ with $\Vert[\o]-\eta\Vert
\leq\e$ is the cohomology class of a Morse form
$\o(\eta)$ such that $v$ is an $\o(\eta)$-gradient.
Choose then the linear forms $\eta_i$ so as to
satisfy Lemma 2.2 and obtain the conclusion.\quad$\square$
\enddemo

\proclaim{Lemma 2.11}
Let $\g,\g'$ be two Morse forms with
$[\g]\not=0, [\g']\not=0$.
 Assume that a vector field
$v$ satisfying the transversality assumption
is  $\g$-gradient and $\g'$-gradient.
Let $x,y\in S(\g), \ind x=\ind y+1$.

Then
\roster\item
If $\sq {\g}=\a [\g']$ with $\a<0$, then
$I(x,y;v)$ is finite.
\item If there is no $\a<0$ with $\sq \g=\a\sq {\g'}$,
then there is an integral ($\sq\g,\sq {\g'})$-cone
$\D$ and $b\in\ZZZ^m$ such that
$Q(I(x,y;v))\subset \D+b$.
\endroster
\endproclaim
\demo{Proof}
1) Obvious.
2) If $\sq \g$ and $[\g']$ are linearly dependant, then
$\sq \g=\a\sq {\g'}$ with $\a>0$ and our lemma follows
from 2.10. Therefore we can assume that $\sq \g$ and $[\g']$
are linearly independant, which imply that there is
$h\in\ZZZ^m$ with $[\g](h),[\g'](h)<0$.
From 2.10 we know that there are integral $\sq \g$-cone
$\G_1$ and $a_1\in\ZZZ^m$ such that
$Q(I(x,y;v))\subset \G_1+a_1$.
Also there are integral $[\g']$-cone $\G_2$
and $a_2\in\ZZZ^m$ such that
such that $Q(I(x,y;v))\subset \G_2+a_2$.
Adding to the generators of $\G_1$ and $\G_2$
some other integral vectors we can assume that
$h\in\Int \G_1,~~h\in\Int \G_2$. Then there exists
$N\in\NNN$ such that
$\G_1+a_1\subset\G_1-Nh$ and $\G_2-Nh\supset \G_2+a_2$.
Thus $Q(I(x,y;v))\subset
(\G_1-Nh)\cap(\G_2-Nh)=\G_1\cap\G_2-Nh$.
The set $\G_1\cap\G_2$ is a $([\g],[\g'])$-cone, and by
2.1 there is an integral $([\g],[\g'])$-cone $\D$ such that
$\G_1\cap \G_2\subset\D$ \th
$Q(I(x,y;v))\subset \D-Nh$.
\quad$\square$\enddemo

By the remark 2.7 the incidence coefficient
$\bbar n(\bbar x,\bbar y;v)\in (\ZZZ[\ZZZ^m])^-_{[\o]}$
is defined. We shall assume that the liftings $\bbar x$
of points $x\in S(\o)$ are chosen so that $Q(\wi x)=\bbar x$.
Note that obviously $\supp(\bbar n(\bbar x,\bbar y;v))
\subset Q(\supp \wi n(\wi x,\wi y;v))$.

\subhead{Proof of Theorem B}
\endsubhead
For a Morse form $\xi$ such that $[\xi]\in H^1(M,\QQQ), [\xi]\not=0$
we denote by $\xi_0$ the (unique)
Morse form, such that $[\xi_0]$ is an indivisible class
in $H^1(M,\ZZZ)$ and that $\xi_0=\mu\xi$ with $\mu>0$.
The map $M\to S^1$, corresponding to $\xi_0$, will be denoted
by $f_0\langle\xi\rangle$.

Define now $\GG t_1(\o)$ as the set of all
$\o$-gradients $v\in\GG t(\o)$, satisfying the following
property:
\vskip0.2in
\it
$({\Cal C})$: There is $\e:0<\e\leq\Vert[\o]\Vert/2$
such that $v$ is an $\O_\e$-gradient and there is a Morse form
$\xi\in\O_\e$ with $[\xi]\in H^1(M,\QQQ)$ such that
$v\in\GG t_1(f_0\langle\xi\rangle)$.
\rm
\vskip0.2in
We shall now prove the properties of $\GG t_1(\o)$.

1) $\GG t_1(\o)$ is $C^0$-open in $\GG t(\o)$.

Indeed, if $v$ satisfies $({\Cal C})$ then every $\o$-gradient
$u$, sufficiently close to $v$, is also an $\O_\e$-gradient
(by 2.5) and $u\in\GG t_1(f_0\langle\xi\rangle)$
since $\GG t_1(f_0\langle\xi\rangle )$ is $C^0$-open
in $\GG t(f_0\langle\xi\rangle)$.

2) $\GG t_1(\o)$ is $C^0$-dense in $\GG t(\o)$.

Indeed, if $v\in\GG t(\o)$, then there is an $\e>0$ such that
$v$ is an $\O_\e$-gradient. Choose any form
$\o'\in\O_\e$ with $[\o']\in H^1(M,\QQQ)$.
 Then by Theorem 1.3
arbitrarily close to $v$ there is
an $\o'$-gradient $u\in\GG t_1(f_0\langle\o'\rangle)$.
By 2.5 $u$ is also an $\O_\e$-gradient.

3) If $v\in\GG t_1(\o)$, then
$\bbar n(\bbar x,\bbar y;v)\in S_{[\o]}^{-1}\zzz$.

Indeed, $v$ is an $\o$-gradient and a $\xi$-gradient for some
$\xi\in\O_\e$ with $[\xi]\in H^1(M,\QQQ)$.
Note that if $[\xi]$ and $[\o]$ are linearly dependant, then
$[\xi]=\a[\o]$ with $\a>0$. (Indeed, $[\xi]=\a[\o]$ with
$\a<0$ would imply
$(1-\a)[\o]+\sum\e_ia_i^*=0$, which contradicts
$\e\leq\Vert[\o]\Vert/2$. \Th 2.11 implies
that there is an integral
$([\xi],[\o])$-cone $\D$, such that
$\supp \bbar n(\bbar x,\bbar y;v)\subset \D+b$
for some $b\in\ZZZ^m$.
Consider $\bbar n(\bbar x,\bbar y;v)$ as an element of
$(\zzz)^-_{[\xi]}$.
Lemma 2.3 together with  Theorem 1.3 imply
 that it belongs to the localization
$S_{[\xi]}^{-1}\zzz$, \th there are $P,Q\in\zzz$
such that $P=Q\cdot\bbar n(\bbar x,\bbar y;v)$.
Since $\bbar n(\bbar x,\bbar y;v)\in \s_\D^{-1}\ZZZ[\D]\sphat~$
the faithful flatness
property imply
 $\bbar n(\bbar x,\bbar y;v)\in\s_\D^{-1}S_\D^{-1}
\ZZZ[\D]\subset S_{[\o]}^{-1}\zzz$.\quad$\square$

\subhead
D. Incidence coefficients with respect to
the universal cover: the exponential estimate
\endsubhead

Passing to the exponential estimate, we need first a lemma.

\proclaim{Lemma 2.12}
Let $\xi,\eta:\RRR^m\to\RRR$ be non zero linear forms, and let
$\G$ be a $(\xi,\eta)$-integral cone. Then there is
$A>0$ such that for every $b\in\RRR^m$ there is $B$ such that
for every $c\in\RRR$ we have:
$$(\G+b)\cap\xi^{-1}([c,\infty[)\subset
(\G+b)\cap\eta^{-1}([Ac+B,\infty[)$$
\endproclaim
\demo{Proof}
Abbreviate $\xi^{-1}([c,\infty[)$ by $\{\xi\geq c\}$. It
is sufficient to prove that there is $A\geq 0$ such that
$$\G\cap\{\xi\geq c\}\subset\G\cap\{\eta\geq Ac\}\tag*$$
(the case of general $b\in\RRR^m$ follows
then with $B=\eta(b)-A\xi(b)$).
To prove (*) let $e_i$ be the generators of $\G$, and choose
$A>0$ such that $\eta(e_i)\geq A\xi(e_i)$ for all $i$.
Then $x\in\G\cap\{\xi\geq c\}$ means: $x=
\sum\l_i e_i$ where $\l_i\geq 0$ and 
$\sum\l_i\xi(e_i)\geq c$.
This implies  $\eta(x)=\sum\l_i\eta(e_i)\geq
\sum\l_iA\xi(e_i)\geq A\xi(x)$.
$\square$\enddemo

Let $\l=\sum_g n_gg\in\ZZZ[\pi_1M]_\xi^-$
and $c\in\RRR$.
Denote by $\l[c]$ the element
$\sum\limits_{\xi(g)\geq c} n_g\cdot g$
of
$\ZZZ[\pi_1M]$, and by $N_c(\l)$ the norm of $\l[c]$, that is
$N_c(\l)=\sum\limits_{\xi(g)\geq c}\vert n_g\vert$.

\subhead{Proof of Theorem C}
\endsubhead
By the definition of $\GG t_1(\o)$ there is
$\e\in]0,\frac 12\Vert[\o]\Vert[$, and
a Morse form
$\xi\in\O_\e$, such that
$[\xi]\in H^1(M,\QQQ)$ and $v\in\GG t_1(f_0\langle\xi\rangle)$.
 \Th there is
an integral $([\xi],[\o])$-cone $\D$ and $b\in\ZZZ^m$,
such that $Q(\supp n(\wi x,\wi y;v))\subset\D+b$, and,
\th , 2.12 implies that the set $\supp \wi n(\wi x,\wi y;v)\cap
(\{\o\})^{-1}([c,\infty[)\subset
\supp \wi n(\wi x,\wi y;v)\cap
(\{\xi\})^{-1}([Ac+B,\infty[)$, and since
$v$ is a $\xi$-gradient, our estimate follows from
Theorem 1.3 and Proposition 1.2. \quad$\square$

\subhead{E. On the Novikov complex for a Morse form} \endsubhead
In this paper we do not use the notion of Novikov complex, working only
with the incidence coefficients. The latter were introduced however
in [5] as the matrix entries of the boundary operators
in the Novikov complex. In this Subsection we use the results
of Subsection B and Subsection C to give a simple proof
of the fact that $\partial^2=0$ in this complex (reducing it to
the corresponding statement about rational Morse forms, proved in
[7]).

We assume here the terminology of Subsection B of the Introduction.
Moreover, if $\D$ is an integral $[\o]$-cone, we denote by
$\L_\D$ the subset of
$\ZZZ[\pi_1M]^-_{\{\o\}}$,
defined by
$$
\L_\D=\{\l\vert \exists b\in\ZZZ^m~\text{\rm such~that~}
Q(\supp \l)\subset\D+b\}
$$
It is not difficult to see that $\L_\D$ is a subring of
$\ZZZ[\pi_1M]^-_{\{\o\}}$.

Now let $v$ be an $\o$-\gr satisfying transversality assumption
and $\e>0$ so small that $\O_\e$ is a Morse family and $v$ is
an $\O_\e$-\gr. Let $\xi\in\O_\e$ be a Morse form such that
$[\xi]\in H^1(M,\QQQ)$. 
Then it follows from 2.11 that there is
an integral $([\xi],[\o])$-cone $\D$, such that for every
$x,y\in S(\o)$ with $\ind x=\ind y+1$ we have:
$ \wi n(\wi x,\wi y;v) \subset  \L_\D$.
Therefore the \ho $\partial : C_*(\o,v)\to C_{*-1}(\o,v)$
is defined actually over the ring $\L_\D$
and to verify that $\partial^2=0$
it is sufficient to verify it over the ring $\ZZZ[\pi_1M]^-_{\{\xi\}}$,
which is done in [7].

\newpage

\head
\S 3. An example
\endhead

In this section we shall construct a Morse map
$M\to S^1$ on a closed 3-\ma $M$, having two critical
points: $x$ of index 2 and $y$ of index 1,
such that $n(\bar x,\bar y;v)=\sum_{k\geq 0}n_kt^k$
with $n_k\sim \a\cdot\b^k$ with $\a<0,\b>0$.
For any $f$-gradient $w$ sufficiently close to $v$
in $C^0$ topology, we shall have:
$n_k(\bar x,\bar y;w)=n_k(\bar x,\bar y;v)$.

The construction is illustrated on the fig.1
and we invite the reader to consult it.
\footnote{The author does not have a TEX-version
of the figure. The figure is available on request
from the author.}
We start with a torus $T^2$, move away from it
two open discs and obtain a surface $S$ with
two components of boundary. Choose and fix a parallel $\b$ and a meridian
$\a$ of this twice punctured torus.
The copies of this surface (resp. the copies of $\a$, $\b$ etc.)
will be denoted by the same letter $S$ (resp. $\a$, $\b$, etc.)
adding indices in order to distinguish between them.
The corresponding discs will be denoted by $D_1,D_2$.
We glue a copy of $S$, denoted by $S(1,1)$ to a copy
of $S$, denoted by $S(1,2)$ and close the boundary by two discs
$D_1(1), D_1(2)$ (as it is shown on the figure 1).
The resulting surface is called $N$.

We glue three copies of $S$ sucsessively, close the boundary
and obtain a closed surface
$L=D_1(1/2)\cup S(1/2,1)\cup S(1/2,0)\cup S(1/2,2)\cup D_2(1/2)$.

The surface $K=D_1(0)\cup S(0,1)\cup S(0,2)\cup D_2(0)$
is constructed similarly to $N$.

Let $W_1$ be the cobordism between $N$ and $L$, corresponding
to the surgery on the circle $\b(1/2,2)$.
Introduce on $W_1$ the corresponding Morse function
$F_1:W_1\to[1/2,1]$ with one critical point
$x$ of index 2,
 and such that
$F_1^{-1}(1)=N, F_1^{-1}(1/2)=L$. We can find an $F_1$-gradient
$v_1$ such that
\roster\item
$\stind {(-v_1)}{1/2}1 $ restricted to
$D_1(1/2)\cup S(1/2,1)\cup S(1/2,0)$
is
a diffeomorphism of this surface onto $D_1(1)\cup S(1,1)\cup S(1,2)$
which identifies $S(1/2,0)$ with $S(1,2)$ and $S(1/2,1)$ with
$S(1,1)$, slightly diminished from the left, so that the image of
$D_1(1/2)$ contains $D_1(1)$ in its interior (and therefore
$\stind {(-v_1)}{1/2}1 (D_1(1))\subset \Int D_1(1/2)$).
\item
$\stind {(-v_1)}{1/2}1 (S(1/2,2)\cup D_2(1/2))=D_2(1)$;
\quad
$\stind {(-v_1)}{1/2}1 (D_2(1/2))\subset\Int D_2(1)$;
and for some $\d>0$ the
$\stind {(-v_1)}{1/2}1$-image of a $\d$-tubular neighborhood
of $\b(1/2,2)$
is in $\Int D_2(1)$ (this image equals to $D_\d(-v_1)\cap F_1^{-1}(1)$).
\endroster

Let $W_0$ be the cobordism between $K$ and $L$, corresponding
to the surgery on the circle $\b(1/2,1)$.
Introduce on $W_0$ the corresponding Morse function
$F_0:W_0\to[0,1/2]$ with one critical point
$y$ of index 1,
 and such that
$F_0^{-1}(1/2)=L, F_0^{-1}(0)=K$.
 We can find an $F_0$-gradient
$v_0$ such that
\roster\item
$\stind {(v_0)}{1/2}0 $ restricted to
$S(1/2,0)\cup S(1/2,2)\cup D_2(1/2)$
is a diffeomorphpism of this surface onto
$D_1(0)\cup S(0,1)\cup S(0,2)$
which identifies $S(1/2,0)$ with $S(0,1)$ and $S(1/2,2)$ with
$S(0,2)$ slightly diminished from the right so that the
image of $D_2(1/2)$ contains $D_2(0)$ in its interior (and \th
$\stind {(-v_0)}{1/2}0 (D_2(0))
\mathbreak
\subset \Int D_2(1/2)$).
\item
$\stind {(v_0)}{1/2}0 (S(1/2,1)\cup D_1(1/2))=D_1(0)$;
\quad
$\stind {(v_0)}{1/2}0 (D_1(1/2))\subset\Int D_1(0)$;
and for some $\d>0$ the
$\stind {(v_0)}{1/2}0$-image of a $\d$-tubular neighborhood
of $\b(1/2,2)$
(this image equals to $D_\d(v_0)\cap F_0^{-1}(0)$)
is in $\Int D_1(0)$.
\endroster

Glue together $W_0$ and $W_1$ along $L$,
denote the resulting cobordism by $W$; we have $\partial W=
K\cup N$. Glue the functions $F_1$ and $F_0$ to a \Mf $F:W\to[0,1]$
and the vector fields $v_1,v_0$ to an $F$-\gr $v$.
We shall now define a ranging system
for $(F,v)$. Consider the set
$\L=\{0,1/2,1\}$ of regular values, and set
$$
\gather A_1=D_1(1);\quad B_1=D_2(1);\\
A_{1/2}=D_1(1/2);\quad B_{1/2}=D_2(1/2);\\
A_0=D_1(0);\quad B_0=D_2(0);
\endgather$$
The fact that $\ran$ is a ranging system for
$(F,v)$ follows immediately from the properties of $v_0,v_1$
cited above.

It is not difficult to compute the \ho
$H(v):H_1(N\sm D_2(1),D_1(1))\to H_1(K\sm D_2(0),D_1(0))$
Namely,
$[\a(1,1)]\mapsto 0,\quad
[\b(1,1)]\mapsto 0,\quad
[\a(1,2)]\mapsto [\a(0,1)],\quad
[\b(1,2)]\mapsto [\b(0,1)].$

It is also obvious that the homology class of
$D(x,v)\cap K$ equals to
$[\b(0,2)]$ and the homology
class of $D(y,-v)\cap N$ equals to
$[\b(1,1)]$.

Consider the embedded curves
$\a(0,1), \b(0,1),\a(0,2),\b(0,2)$
in $K$. Their homology classes form a symplectic basis in
$K$. Consider the isomorphism of $H_1(K,\ZZZ)\approx\ZZZ^4$
given in this basis by the matrix

$$
\left( \matrix 0 & 2 & 1 & 0\\
0 & 0 & 0 & 1\\
0 & 1 & 0 & 0\\
-1 & 3 & 0 & -2\\
\endmatrix
\right).
$$

It preserves the intersection form,
\th it can be realized by a diffeomorphism
$\Phi:K\to K$ and it is easy to see that we can assume
that $\Phi(x)=x$ for
$x\in D_1(0)\cup D_2(0)$. Denote by $\Psi$ the composition
of $\Phi$ with the subsequent identification of $K$ with $N$.

Denote by $M$ the 3-\ma obtained by glueing
of $K$ to $N$ by means of $\Psi$, and let $f:M\to S^1$
be the  \Mf obtained from $F$.
The corresponding cyclic covering $\bar M$ is a union of
countably many copies of
$W$, denoted by $W[i],i\in\ZZZ$, glued together by the diffeomorphisms
$\Psi:K[i]\to N[i-1]$ of the components of boundaries.

The identification $W[k]\to W$ will be denoted by $J_k$. For $X\subset W$
we denote $J_k^{-1}(X)$ by $X[k]$. Define a lifting
$\bar F:\bar M\to\RRR$ of $f$ by setting
$\bar F\vert W[k]= F\circ J_k +k$. Set
$\Sigma=\{n/2~\vert n\in\ZZZ\}$;
for $n\in\ZZZ$ set $A_n=D_1(0)[n], B_n=D_2(0)[n]$;
for $n=k+1/2, k\in\ZZZ$ set $A_n=D_1(1/2)[k], B_n=D_2(1/2)[k]$.
It is obvious that $\rans$ is a $t$-equivariant ranging system
for $(\bar F,v)$. Set $\bar x=J_0^{-1}(x), \bar y=J_0^{-1}(y)$.

We have: $n_0(\bar x,\bar y;v)=0,
 n_1(\bar x,\bar y;v)=\Psi_*([\b(0,2)])~\sharp~
[\b(1,1)]=0$. Denote by $D$ the matrix of $\Psi_*\circ H(v)$; then
$$D=\left(
\matrix
0 & 0 & 0 & 2\\
0 & 0 & 0 & 0\\
0 & 0 & 0 & 1\\
0 & 0 & -1 & 3
\endmatrix
\right).
$$
Then for $k\geq 1$ we have
$n_{k+1}(\bar x,\bar y;v)=(\Psi_*\circ H(v))^k(\Psi_*([\b(0,2)])
~\sharp~ [\b(1,1)]$.
To abbreviate the notation we denote
$[\a(1,1)]$ by $a_1$,
$[\a(1,2)]$ by $a_2$,
$[\b(1,1)]$ by $b_1$,
$[\b(1,2)]$ by $b_2$.
Then $n_{k+1}(\bar x,\bar y;v)$ is the coefficient of
$a_1$ in $D^k(b_1-2b_2)=(-2)D^k(b_2)$.
To find $D^k(b_2)$ assume by induction that
$D^k(b_2)=\a_ka_1+\g_k a_2+\b_kb_2$. Then
$D^{k+1}(b_2)=-\g_kb_2+\b_l(2a_1+3\b_2+a_2)=
2\b_ka_1+(3\b_k-\g_k)b_2+\b_ka_2$. \Th the vectors
$(\b_k,\g_k)$ satisfy
$$
\pmatrix
\b_{k+1}\\
\g_{k+1}
\endpmatrix
=
\pmatrix
3 & -1\\
1 & 0
\endpmatrix
\cdot
\pmatrix
\b_{k}\\
\g_{k}
\endpmatrix
\quad\text{\rm and }\quad
\pmatrix
\b_0\\
\g_0
\endpmatrix
=
\pmatrix
1\\
0
\endpmatrix
.$$
The coefficient $\a_{k+1}=2\b_k$.

Denote
$$\pmatrix
3 & -1\\
1 & 0
\endpmatrix
$$ by $C$; then the explicit computation shows
$$C^n=\frac 1{\sqrt 5}\cdot
\pmatrix
\l_1^{n+1}-\l_2^{n+1} & \l_2^{n}-\l_1^{n}\\
{} & {}\\
\l_1^{n}-\l_2^{n} & \l_2^{n-1}-\l_1^{n-1}
\endpmatrix
$$
where $\l_1=\frac {3+\sqrt 5}2, \quad \l_2=\frac {3-\sqrt 5}2$.
\Th $$n_{k+1}(\bar x,\bar y;v) =
-\frac 4{\sqrt 5}\cdot \bigg( \bigg(\frac {3+\sqrt 5}2\bigg)^k -
 \bigg(\frac {3-\sqrt 5}2\bigg)^k\bigg).$$

\newpage

\head
\S 4. Exponential estimates of the absolute number
of trajectories.
\endhead

\subhead
A. Bunches of ranging systems
\endsubhead

Let $f:W\to[a,b]$ be a \Mf on a \co $W$,
$f^{-1}(b)=V_1, f^{-1}(a)=V_0, \dim W=n$.
Let $v$ be an $f$-\gr, $\L=\{\l_0,...,\l_k\}$
be a set of regular values of $f$, such that
$\l_0=a<\l_1<...<\l_{k-1}<\l_k=b$
and between any two adjacent values $\l_i$ and $\l_{i+1}$
there is exactly one critical value of $f$.

\definition{Definition 4.1}
Assume that for each integer $s:0\leq s\leq n$ and every
$\l\in\L$ there are given compacts
$A_\l^{(s)}, B_\l^{(s)}$ in $f^{-1}(\l)$. We shall say that
$\rran$ is a \it bunch of ranging systems for $(f,v)$ \rm
(abbreviation: BRS for $(f,v)$) if the following conditions hold:
 \roster\item $A_\l^{(0)}= B_\l^{(0)}=\emptyset$
\item $r\leq s\Rightarrow \bigg( A_\l^{(r)}\subset
A_\l^{(s)},  B_\l^{(r)}\subset
B_\l^{(s)}\bigg)  $
\item $ A_\l^{(r)}\cap B_\l^{(s)}=\emptyset$ for $r+s\leq n$
\item If $\l,\m\in\L$ are adjacent, then

$\stind v\mu\lambda (A_\m ^{(r)})\subset \Int A_\lambda^{(r)}
 \quad\text{and}
\quad
\stind {(-v)}\lambda\mu (B_\l^{(s)})\subset
\Int B_\m^{(s)}$ for every $r,s$.

\item  Let $\lambda,\mu\in\Lambda$ be adjacent. Then for
 every $p\in S(f)\cap \fpr f{[\lambda,\mu]}$
we have:

i) $D(p,v)\cap \fpr f{\lambda}\subset \Int A_\lambda^{(\ind p)}$ and

ii)$D(p,-v)\cap \fpr f{\mu}\subset \Int B_\mu^{(n-\ind p)}$.
\quad $\triangle$
\endroster
\enddefinition

Let $\rran$ be a BRS for $(f,v)$. The following properties are
either obvious or
proved similarly to the corresponding properties of
ranging systems (see \S 4 of [9]).

\roster\item"i)" For every $s:0\leq s\leq n$ the family
$\{(A_\l^{(s)}, B_\l^{(n-s-1)})\}_{\l\in\L}$
is a ranging system for $(f,v)$.
\item"ii)" $v$ is an almost good $f$-gradient
\item"iii)" There is $\e>0$ such that for every
$f$-gradient $w$ with $\Vert w-v\Vert<\e$
the system $\rran$ is a BRS for $(f,w)$.
\item"iiii)" For every $\d>0$ sufficiently small
the following strengthening of 5) holds: \endroster
$$\left\{
\gathered
D_\d(p,v)\cap \fpr f{\lambda}\subset \Int A_\lambda^{(\ind p)}
 \aand
\\
D_\d(p,-v)\cap \fpr f{\mu}\subset \Int B_\mu^{(n-\ind p)}
\endgathered\right.
\tag5\d
$$
For $p\in S(f)$ we shall denote by $\l_p$ resp. $\m_p$
the maximal resp. minimal element of $\L$ with the property
$\l<f(p)$(resp. $\m>f(p)$).

Let $\rran$ be a BRS for $(f,v)$. Choose $\rho>0$ so small
that for every $p\in S(f)$ we have:
$\l_p<\l_p+\rho<f(p)-\rho<f(p)+\rho<\m_p-\rho.$
Let $\d>0$ satisfy  (5$\d$)
above and assume further that
$$
\left\{
\gathered
(1)
\text{\rm ~~If~} p,q\in S(f) \text{~\rm and~} f(p)=f(q) \text{~\rm then~}
D_\d(p,v)\cap D_\d(q,v)\cap f^{-1}([\l_p,\m_p])=\emptyset  \\
\aand
D_\d(p,-v)\cap D_\d(q,-v)\cap f^{-1}([\l_p,\m_p])=\emptyset.\\
(2)
\text{\rm ~~For~ every~} p\in S(f) \text{\rm ~we have~}
D_\d(p)\subset f^{-1}(]f(p)-\rho,f(p)+\rho[).
\endgathered
\right.
\tag{\d}
$$
Let $\e>0$ satisfy (iii) above.
The \ma $S_p=D(p,-v)\cap f^{-1}(\m_p)$ is an embedded
sphere of dimension $n-\ind p-1$. The normal bundle to $S_p$ in
$f^{-1}(\m_p)$ is trivial. Chose an embedding
$\psi:S_p\times B^{\ind p}(0,\t)\to \ffmin (\m_p)$
(where $\t>0$) which is a diffeomorphism onto its image,
and such that
$\psi\vert (S_p\times\{0\})=\id$.

For $\k\in]0,\t]$ the set
$\psi(S_p\times B^{\ind p}(0,\k))$
will be denoted by $\Tub(S_p,\k)$ and for $\eta\in B^{\ind p}(0,\t)$
the embedded sphere $\psi(S_p\times \{\eta\})$ by
$S_p(\eta)$.

Assume that $\t$ is so small that
$\overline{\Im \psi}\subset B_\d(p,-v)$. Let $\tau>0$ and denote by
\break
$\Psi:[-\tau,0]\times S_p\times B^{\ind p}(0,\t)\to W$
the map
$(t,s,h)\mapsto\g(\psi(s,h),t;v)$.
It is a diffeomorphism onto its image.
We shall assume that $\tau$ is so small that
$\Im\Psi\subset \ffmin(]\m_p-\rho,\m_p])$.

The vector field $\Psi^{-1}_*(v)$
equals to $(1,0,0)$.

Let $\xi\in B^{\ind p}(0,\t/2)$ and
let $H_t\langle\xi\rangle, t\in[-\tau,0]$ be an isotopy
of $B^{\ind p}(0,\t)$ having the following properties:
\roster\item
 $\frac d{dt} H_t\langle\xi\rangle(x)=0$
for $t\in[-\tau,-{\frac 23}\tau]\cup [-{\frac 13}\tau,0]$.
\item $H_t\langle\xi\rangle(x)=x$ for $x\in B^{\ind p}(0,\t)\sm
B^{\ind p} (0,\t/2)$, and for $t\in [-\tau,-{\frac 23 }\tau]$.
\item $H_0\langle\xi\rangle (0)=\xi$.
\endroster

Then $\frac d{dt} H_t\langle\xi\rangle$ is a time-dependant
vector field on $B^{\ind p}(0,\t)$.

Define a vector field $w\langle\xi\rangle$ on
$[-\tau,0]\times S_p\times B^{\ind p}(0,\t)$
by
$w\langle\xi\rangle(t,x,y)=(1,0,\frac d{dt} H_t\langle\xi\rangle)$.

Note that

(A) $\forall x\in S_p\times B^{\ind p}(0,\t/2)$ we have
$\g(x,-\tau;w\langle\xi\rangle)\in \{-\tau\}\times S_p\times
B^{\ind p}(0,\t/2)$.

Choose $\a\in]0,\t/2[$ so small that for each
$\xi\in B^{\ind p}(0,\a)$ there is an isotopy
$H_t\langle\xi\rangle$ such that
$$\Vert v-\Psi_*(w\langle\xi\rangle)\Vert<\e.\tag*$$

We perform this construction for every
 $p\in S(f)$. We distinguish between
the notation of
the corresponding objects
by adding in the index the notation of the point, e.g.
$\Psi_p,\xi_p,w_p\langle\xi_p\rangle,...$
We assume that $\t_p$, resp. $\a_p,\tau_p$ are chosen to be independant
of $p$, and we denote them simply by $\t,\a,\tau$.

Chose now for every $p\in S(f)$ a vector $\xi_p\in B^{\ind p}(0,\a)$,
and the corresponding vector field
$w_p\langle\xi_p\rangle$, satisfying (*).
The set $\{\xi_p\}_{p\in S(f)}$ will be denoted by $\vec\xi$.

Define a new vector field $u=v\langle\vec\xi\rangle$ setting
$$\left\{
\gathered
u(x)=v(x) \iif x\in W\sm \cup_p\Im\Psi_p \\
u(x)=(\Psi_p)_*(w_p\langle\xi_p\rangle)(x) \ffor x\in\Im\Psi_p
\endgathered
\right.
$$

It follows from the construction that for every $\vec\xi$
the vector field $v\langle\vec\xi\rangle$
is an $f$-\gr.
Note also that for $p\in S(f)$ we have
$D(p,-v\langle\vec\xi\rangle)\cap \ffmin(\m_p)=
S_p(\xi_p)$.
The condition (*) implies that $\rran$ is
a BRS for $(f,v\langle\vec\xi\rangle)$; therefore
$v\langle\vec\xi\rangle$ is an almost good $f$-\gr.
It will be called \it regular perturbation
of $v$ corresponding to $\vec\xi$. \rm

\proclaim{Proposition 4.2}
\roster\item Let $p,q\in S(f),\ind p=\ind q+1$. The set
of $(-v\langle\vec\xi\rangle)$-trajectories
joining $p$ with $q$\footnote{We identify here two trajectories
which differ by a parameter change.}
is in a bijective correspondence with
$(D(p,v)\cap\ffmin(\m_q))\cap S_q(\xi_q)$.
\item
The vector field $v\langle\vec\xi\rangle$ is a good $f$-\gr
if for every
$ p,q\in S(f)$ with $\ind p=\ind q+1$
the \sma $D(p,v)\cap\ffmin(\m_q)$ of $\ffmin(\m_q)$ is transversal to
$S_q(\xi_q)$.
\endroster
\endproclaim
\demo{Proof}
1) Let $\g$ be a $(-v\langle\vec\xi\rangle)$-trajectory, joining
$p$ with $q$. Let
$\g(t_0)\in\ffmin(\m_q)$, where $t_0\in\RRR$. I claim that for
$t<t_0$ we have $\g(t)\notin\supp (v-\vxi)$.
Indeed, if the opposite is true, let $t_1<t_0$ be the first moment when
$\g$ intersect $\supp(v-\vxi)$. Then there is
$s\in S(f)$ such that
$\g(t_1)\in\Tub (S_s,\t/2)$. Note that
$\Tub(S_s,\t/2)\subset D_\d(s,-v)\cap \ffmin(\m_s)
\subset \Int B_{\m_s}^{(n-\ind s)}$,
\th $\g(t_1)$ is in the intersection of
$A_{\m_s}^{(\ind p)}$ with $B_{\m_s}^{(n-\ind s)}$;
and $\ind p\geq \ind s+1$.

Further, $f(s)>f(q)$ and $\l_s\geq\m_q$.
Denote by $t_2$ the moment when $\g$ intersects
$\ffmin(\l_s)$.
Then $t_0\geq t_2>t_1+\tau$, and it is easy to see that
$\g\vert[t_1+\tau,t_2]$ does not intersect
\break
$\supp(v-\vxi)$. The property (A) implies that
$\g(t_1+\tau)\in D_\d(s,-v)$, and, therefore,
$\g(t_2)\in D_\d(s,v)\cap\ffmin(\l_s)
\subset \Int A_{\l_s}^{(\ind s)}$, which
implies
$\g(t_0)\in\Int A_{\m_q}^{(\ind s)}$. \Th
$A_{\m_q}^{(\ind s)}\cap D_\d(q,-v)\not=\emptyset$ and, further,
$A_{\m_q}^{(\ind s)}\cap B_{\m_q}^{(n-\ind q)}\not=\emptyset$,
which implies
$\ind s>\ind q$ and $\ind p\geq\ind s+1>\ind q+1$, contradiction.

\Th $\g\vert ]-\infty, t_0]$ is a $(-v)$-\tr.
In the same way one can show that every
$(-v)$-\tr joining $p$ with a point of
$D_\d(q,-v)\cap\ffmin (\m_q)$ does not intersect
$\supp(v-\vxi)$. That proves 1).

2) Note that $\vxi$ is almost good.
\Th $\vxi$ is good if and only if
$$\bigg(\ind p = \ind q+1\bigg)
\Rightarrow
\bigg(D(p,\vxi)\pitchfork D(q,-\vxi)\bigg).$$
Let $p,q\in S(f),\ind p = \ind q+1$. It suffices to prove that
$D(p,\vxi)\cap \ffmin(\m_q)$ is transversal to
$S_q(\xi_q)$. Let $x$ be a point in the
intersection of these manifolds.
In the part 1) we have proved that there is a
$(-v)$-\tr, joining $p$ with $x$ and not intersecting
$\supp(v-\vxi)$. Then a small neighborhood of this \tr
does not intersect $\supp(v-\vxi)$
and the transversality sought follows
from
$(D(p,v)\cap\ffmin(\m_q))\pitchfork S_q(\xi_q)$.
$\quad\square$
\enddemo

\subhead
B. Volume estimates
\endsubhead

We assume here the terminology of the previous subsection.
Let $\rran$ be a BRS for $(f,v)$. Assume that $\d>0$ satisfies
$(\d)$ from Subsection A. Fix an integer $s: 0\leq s\leq n$.
Let $\l<\m$ be adjacent elements of $\L$.
The set
$${\Cal D}_\m^{(s)}=
\ffmin(\m)\sm
\bigg(\Int A_\m^{(s)}\cup \Int B_\m^{(n-s-1)}
\cup \big( \bigcup\limits_{p\in S'}
B_\d(p,-v)\cap\ffmin(\m)\big)\bigg)$$
(where $S'$ stands for the subset of all $p\in S(f)$
such that $f(p)\in ]\l,\m[$ and $\ind p\leq s$)
is a compact subset of the domain of definition of
$\stind v\m\l$. Denote by
$N_\m^{(s)}$ the norm of the derivative of
$\stind v\m\l$ restricted to ${\Cal D}_\m^{(s)}$
(that is $N_\m^{(s)}=
  \underset       { x\in {\Cal D}_\m^{(s)}   }   \to \sup
\Vert (\stind v\m\l)'(x)\Vert$)
and by $A$ the maximum of $(N_\m^{(s)})^k$ over all
$\m\in\L, 0\leq s\leq n$, and $0\leq k\leq n$.

For $p\in S(f)$ denote by $B_{(p)}$ the norm of the derivative of
$\psi^{-1}_p\vert \psi_p(S_p\times D^{\ind p}(0,\t/2))$
(where the \ma $S_p\times B^{\ind p}(0,\t)$
is endowed with the product riemannian metric). Denote by $B$
the maximum of $B_{(p)}$ over all $p\in S(f)$.

\proclaim{Lemma 4.3}
Let $\l,\m\in\L,\l<\m$. Let $s$ be an integer,
$0\leq s\leq n$.
Let $N$ be a \sma of
$\ffmin (\m)\sm B_\m^{(n-s-1)}$, such that
$N\sm A_\m^{(s)}$ is compact. Denote by $k$ the number of elements
of $\L$ belonging to $[\l,\m]$.
Then $N'=\stind v\m\l (N)$
is a \sma of
$\ffmin (\l)\sm B_\l^{(n-s-1)}$, such that
$N'\sm A_\l^{(s)}$ is compact, and

1) $\vvol (N'\sm A_\l^{(s)})\leq A^{k-1}\cdot \vvol (N\sm A_\m^{(s)});   $

2) If $p\in S(f)$ with $\ind p=s$ and  $\m_p=\l$, then
$$\vvol (\psi^{-1}_p(N'\cap\text{\rm Tub}(S_p,\t/2))\leq
BA^{k-1}\vvol(N\sm A_\m^{(s)})$$

\endproclaim
\demo{Proof}
For the proof of 1) note that it suffices to consider the case
$k=2$. For this case we have obviously
$N'\sm A_\l^{(s)}\subset \st v (N\cap {\Cal D}_\m^{(s)})$.
2) follows from 1) since $\text{\rm Tub}(S_p,\t/2)\subset B_\d(p,-v)$
does not
intersect $A_\l^{(s)}$.
$\quad\square$
\enddemo

\subhead
C. Bunches of equivariant ranging systems
\endsubhead

We assume here the terminology of Subsection A of
Introduction. Further, we assume that the homotopy
class of $f$ is indivisible that is $F(xt)=F(x)-1$.

\definition{Definition 4.4}
Let $\Sigma$ be the set of regular values of $F$
such that
\roster\item"$(\Sigma 1)$" $\s\in\Sigma\Rightarrow
\s+n\in\Sigma$ for every $n\in\ZZZ$.
\item"$(\Sigma2)$" If $\s_1,\s_2\in\Sigma$ are adjacent,
there is exactly one critical value of $F$
in $]\s_1,\s_2[$.
\item"$(\Sigma 3)$" For every $A,B\in\RRR$ the set
$\Sigma\cap [A,B]$ is finite.           \endroster

Assume that for every integer $s:0\leq s\leq n$
and every $\s\in\Sigma$ there are given compacts
$A_\s^{(s)}, B_\s^{(s)}$ in $F^{-1} (\s)$. We shall say that
$\rrans$ is a
\it
bunch of $t$-equivariant ranging systems
\rm
(abbreviation: BERS for $(F,v)$)
if:
\roster\item
For every $\mu,\nu\in\Sigma,\mu<\nu$ the system
$\{(A_\s^{(s)} , B_\sigma^{(s)})\}_{\s\in\Sigma,
\mu\leqslant\sigma\leqslant\nu}$
 is a BRS for
\break
$(F\vert\fpr F{[\mu,\nu]},v)$.
\item
$
A_{\sigma -n}^{(s)} = A_\sigma^{(s)}\cdot t^n,\quad
B_{\sigma -n}^{(s)} = B_\sigma^{(s)}\cdot t^n$ ~for every $n\in \ZZZ$.
\quad$\triangle$              \endroster
\enddefinition

Up to the end of this subsection we assume that
$(F,v)$ has a BERS
\break
 $\rrans$.
Choose any $\s\in\Sigma$, denote by $W$ the \co
$\fmin ([\s,\s+1])$; denote $\Sigma\cap [\s,\s+1]$ by
$\L$.
We apply the constructions of Subsections A,B to $W$ and then extend
the results to the whole of $\bar M$ in the $t$-invariant way,
thus obtaining the sets 
$\text{\rm Tub}(S_q,\varkappa)\subset F^{-1}(\m_q)$
for every $q\in S(F)$.
 Let $p,q\in S(F)$; assume that
$F(p)>F(q)$ and $\ind p=\ind q+1$. Denote by $N_{p,q}$
the \sma of dimension $\ind q$
of $S_q\times B^{\ind q}(0,\t/2)$, defined by
$$N_{p,q}=\psi^{-1}_q\bigg( \big(D(p,v)\cap F^{-1}(\m_q)\big)
\cap \big(\text{\rm Tub}(S_q,\t/2)\big)\bigg).$$

The next lemma follows from 4.3.

\proclaim{Lemma 4.5}
There are constants $C,D>0$ such that for every
$p,q\in S(F)$ with
$\ind p=\ind q+1$ and $F(p)>F(q)$
we have:
$\vvol(N_{p,q})\leq C\cdot D^{[F(p)-F(q)]}.$
\quad$\square$
\endproclaim
The next lemma is a direct consequence of Sard Theorem.

\proclaim{Lemma 4.6}
Let $q\in S(F)$. Then there is a residual subset
$Q\subset B^{\ind q}(0,\a)$ such that for every
$\xi\in Q$ and every $p\in S(F)$ with $\ind p=\ind q+1$
we have:
$D(p,v)\pitchfork S_q(\xi)$.
\quad$\square$     \endproclaim

The next proposition follows from 4.5 by the argument
due to V.I.Arnold (see [1, p.81]).

\proclaim{Proposition 4.7}
Let $q\in S(F).$   Then there is a subset
$Q\subset B^{\ind q}(0,\a)$ of full measure such that for
every $\xi\in Q$ we have:

For every $p\in S(F)$ with $\ind p=\ind q+1$ there are constants
$K,R>0$, such that
for every integer $l\geq 0$ we have:
$$\sharp(N_{pt^{-l},q}
\cap (S_q\times\{\xi\}))\leq K\cdot R^l.\quad\square$$

\endproclaim

\proclaim{Corollary 4.8}
Let $\nu>0$. Then there is a good $f$-\gr $u$
with $\Vert u-v\Vert<\nu$ and  constants
$K,L>0$ such that for every $p,q\in S(f)$ with $\ind p=\ind q+1$
we have
$N_l(p,q;v)\leq L\cdot K^l$.
\endproclaim
\demo{Proof}
The construction of regular perturbation, described in
Subsection A, applied to $v\vert W$, gives for every
$\vec\xi=\{\xi_p\}_{p\in S(F)\cap W}, \xi_p\in B^{\ind p}(0,\a)$
an $F\vert W$-\gr $\vxi$. Since
$\supp (v-\vxi)$ does not intersect $\partial W$,
we can extend $\vxi$ $t$-invariantly to $\bar M$ and obtain a $t$-invariant
$F$-\gr, which will be denoted by the same symbol $\vxi$.
Note that since $\rrans$ is a BERS for $(F,\vxi)$,
the vector field $\vxi$ is an almost good $F$-\gr.

If all the $\xi_q\in B^{\ind q}(0,\a)$
are sufficiently small, we shall have
$\Vert u-\vxi\Vert<\nu$. Propositions 4.6, 4.7 imply that
we can choose the components $\xi_q$ of $\vec\xi$
in such a way that
\roster\item For every $p\in S(F),\ind p=\ind q+1$ we have:
$D(p,v)\cap F^{-1}(\m_q)\pitchfork S_q(\xi_q)$.
\item There are $K,R>0$ such that for every
$p\in S(F)$ with $\ind p=\ind q+1$ we have:
$\sharp(D(pt^{-l},v)\cap S_q(\xi_q))\leq K\cdot R^l$.   \endroster

The proposition 4.2 implies then that $u=\vxi$ satisfy
the conclusion. \quad$\square$
\enddemo

\subhead
D. Proof of Theorem D
\endsubhead
We can assume that $S(f)\not=0$ and that the homotopy
class of $f$ is indivisible.
In view of Corollary 4.8 it is sufficient to prove that
the set of $f$-\grs, having a BERS, is $C^0$ dense
in $\GG (f)$. To prove that,
let $v$ be any $f$-\gr and $\e>0$.
Choose any regular value $\l$ of $F$ and denote
by $W$ the cobordism $\fmin ([\l,\l+1])$.
Choose any set of regular values of $F$,
satisfying $(\Sigma 1) - (\Sigma 3)$ of Def. 4.4.
Theorem 4.13 of [9] imply that there is an
almost good $F\vert W$-gradient $u$ and a
ranging pair $(\VV_0, \VV_1)$ for
$(F\vert W,u)$ such that
$\Vert v-u\Vert\leq \e$
and $v=u$ near $\partial W$.
We can also assume that $\VV_0t=\VV_1$.
 Then the procedure, described in
Construction 4.10 of [9]
define a bunch of ranging systems for
$(F\vert W,u)$. Extend $u$ to a $t$-invariant $F$-\gr.
The BRS constructed is easily extended to a BERS
for $(F,u)$. \quad$\square$

We mention here an obvious corollary of
Theorem D which will be of use in the proof of
Theorem E.

For a good $f$-\gr $v, c\in \RRR$ and $x,y\in S(f)$ with
$\ind x=\ind y+1$ we denote by $N_{\geq c}(x,y;v)$
the sum $\sum\limits_{-k\geq c} N_k(x,y;v)$.

\proclaim{Corollary 4.9}
Let $v\in \GG_0(f)$. Then there are constants $R,Q>0$ such that
$\forall x,y\in S(f)$ with $\ind x=\ind y+1$ and
$\forall c\in\RRR$ we have:
$N_{\geq c}(x,y;v)\leq R\cdot Q^{-c}$.\quad$\square$ \endproclaim

\subhead
E. Proof of the Theorem E
\endsubhead

By Lemma 2.5 find $\e>0$ and $\d>0$ such that
$\O_\e$ is a Morse family, $v$ is an $\O_\e$-\gr and
 every $\o$-\gr $u$ with $\Vert u-v\Vert<\d$
 is an $\Omega_\e$-\gr.
Choose some $\xi\in\Omega_\e$ with
$[\xi]\in H^1(M,\QQQ)$ and choose (by Theorem D) a good $\xi$-\gr $u$,
satisfying 2) of Theorem D and such that $\Vert u-v\Vert<\d$;
$u$ is then an $\O_\e$-\gr.
 By Lemma 2.11
for every $x,y\in S(\om)$ there is an integral
$([\xi],[\om])$-cone
$\D$ and $b\in \ZZZ^m$ such that
$Q(I(\tilde x,\tilde y;u))\subset \D+b$.
Lemma 2.12 imply that there are $A,B$ such that
for every $\l\in\RRR$ we have
$I(\tilde x,\tilde y;u)\cap\{\om\}^{-1}([\l,\infty[)
\subset
I(\tilde x,\tilde y;u)\cap\{\xi\}^{-1}([A\l+B,\infty[)$,
and this together with 4.9 implies the conclusion.
\quad$\square$

\enddocument